\documentclass[reqno,10pt]{amsart}
\makeatletter
 \usepackage{enumitem}

 \usepackage{amsmath}
\allowdisplaybreaks[4]

\newcommand{\Rmnum}[1]{\expandafter\@slowromancap\romannumeral#1@}

\newtheorem{lemma}{Lemma}[section]

\newtheorem{theorem}{Theorem}[section]
\newtheorem{corollary}{Corollary}[section]

\usepackage{color}
\usepackage{cite}
\usepackage{mathrsfs}
\numberwithin{equation}{section}

\title[Compressible Euler equations]
{The time asymptotic expansion for the compressible Euler equations with damping}%
\author[F. Huang and X. Wu]{}

\email{fhuang@amt.ac.cn}
\email{xcwu22@csu.edu.cn}
\thanks{* Corresponding author.}

\begin{document}

\maketitle \centerline{\scshape    Feimin Huang$^{a,b}$ \ \    Xiaochun Wu$^{*,c}$}
\medskip
{\footnotesize
 \centerline{$^a\ \!$Academy of Mathematics and systems Science, Chinese Academy of Sciences}
   \centerline{Beijing 100190, China}
    \centerline{$^b\ \!$School of Mathematical Sciences, University of Chinese Academy of Sciences}
 \centerline{Beijing 100049, China}
    \centerline{$^c\ \!$School of Mathematics and Statistics, HNP-LAMA, Central South University }
    \centerline{Changsha 410083, China}}

\medskip
\bigskip
\begin{abstract}
In 1992, Hsiao and Liu \cite{Hsiao-Liu-1} firstly showed that the solution to the compressible Euler equations with damping  time-asymptotically converges to the diffusion wave $(\bar v, \bar u)$ of the porous media equation. In \cite{Geng-Huang-Jin-Wu}, we proposed a time-asymptotic expansion around the diffusion wave $(\bar v, \bar u)$, which is a better asymptotic profile than $(\bar v, \bar u)$.
In this paper, we rigorously justify the time-asymptotic expansion by the approximate Green function method and the energy estimates. Moreover, the large time behavior of the solution to compressible Euler equations with damping is accurately characterized by the time asymptotic expansion.
 \end{abstract}

\section{Introduction}
In this paper, we  are concerned with the  1-d compressible
Euler equations with damping in Lagrangian coordinates, which reads as
\begin{align}\label{bE:2.1}
\begin{cases}
v_t-u_{x}=0,\\
u_t+P(v)_{x}=-\alpha u,
\end{cases}
\end{align}
with the initial data
\begin{align}\label{E:1.2}
(v,u)(x,0)=(v_0,u_0)(x)\rightarrow (v_{\pm},u_{\pm}) \quad as\quad x\rightarrow \pm \infty, 
\end{align}
where $v=v(x,t)$ denotes the specific volume,  $u=u(x,t)$ is the velocity and $P(v)$ is
the pressure satisfying $P(v)>0, P^{\prime} (v)<0$. The damping term $\alpha u$ is the friction effect with the physical parameter $\alpha>0$.  Without loss of generality, we assume $\alpha=1$ in this paper.

Since the velocity $u_t$ in $\eqref{bE:2.1}_2$ decays to zero faster than the other terms due to the damping effect, it is commonly conjectured that \eqref{bE:2.1} is time-asymptotically equivalent to the porous medium equation (PME):
\begin{align}\label{E:1.3}
\begin{cases}
 \bar{v}_t-\bar u_{x}=0, \\
P(\bar v)_{x}=-\bar{u},\quad \mbox{Darcy's law},
\end{cases}
\quad \mbox{equivalently,}\quad 
 \bar{v}_t=-P(\bar v)_{xx},\quad \mbox{PME},
\end{align}
and the velocity $\bar u$ obeys the Darcy's law.
This conjecture was firstly justified by Hsiao and Liu \cite{Hsiao-Liu-1} around the diffusion wave $\bar v(x,t)=:\bar{v}(\xi), \xi=\frac{x}{\sqrt{1+t}}$  of PME with the boundary condition at the far field
\begin{align}\label{E:2.1new}
\bar v(x,t)\rightarrow v_{\pm} \hspace{0.2cm}as\quad x \rightarrow \pm \infty.
\end{align}
Since $\bar v(\xi)$ is the self-similar solution to  the PME \eqref{E:1.3}, then $\bar v(\xi)$ solves the following ODE:
 \begin{align}\label{a1}
 \frac{1}{2}\xi \bar v^{\prime}-(P({\bar v}))^{\prime\prime}=0, \quad ^{\prime}=\frac{d}{d\xi}
\end{align}
with
 \begin{align}
 \lim_{\xi\rightarrow\pm\infty} {\bar v}(\xi)=v_{\pm}.
 \end{align}
 In fact, there exists a unique solution $\bar v(\xi)$ up to a shift in \cite{Hsiao-Liu-1}.
Once $\bar v(\xi)$ is determined, $\bar u(x,t)$ is defined as 
\begin{equation}\label{a2}
\bar u(x,t)=-(1+t)^{-\frac{1}{2}}(P(\bar v))^{\prime}.
\end{equation}
Furthermore, it was shown in \cite{Hsiao-Liu-1} that the solution $(v, u)(x,t)$ to \eqref{bE:2.1}-\eqref{E:1.2} converges to the diffusion wave $(\bar v, \bar u)(x,t)$ given in \eqref{a1}-\eqref{a2} in the form of $$\| v-\bar v,  u-\bar u\|_{L^{\infty}}=O(1)((1+t)^{-\frac12},(1+t)^{-\frac12}).$$ Later,  Nishihara \cite{Nishihara} improved the decay rate to $$\| v-\bar v,  u-\bar u\|_{L^{\infty}}=O(1)((1+t)^{-\frac34},(1+t)^{-\frac54}).$$ By constructing a fine approximate Green function and elaborate estimates, Nishihara-Wang-Yang \cite{Nishihara-Wang-Yang} 
 further improved  to $$\| v-\bar v,  u-\bar u\|_{L^{\infty}}=O(1)((1+t)^{-1},(1+t)^{-\frac32}).$$ Mei \cite{Mei} 
 chose another asymptotic profile (not self-similar solution $\bar v(\frac{x}{\sqrt{1+t}})$ of PME) and obtained the corresponding decay rate.
For the other interesting works, see \cite{ Hsiao-Liu-2,Hsiao-Luo,  Zhao,Zheng,Lions-Perthame-Tadmor, Lions-Perthame-Souganidis, Liu-Yang-1997, Liu-Yang-2000,Geng-Huang, Huang-Marcati-Pan, Huang-Pan-03,Huang-Pan-06, Huang-Pan-Wang, Liu, Luo-Zeng, Serre-Xiao, Sugiyama-1,Zhu-03} and the references therein. 

On the other hand, we in \cite{Geng-Huang-Jin-Wu} considered the compressible Euler equation with time-dependent damping
\begin{align}\label{bE:2.1nb}
\begin{cases}
v_t- u_{x}=0,\\
u_{t}+P( v)_{x}=-\frac{1}{(1+t)^{\lambda}}  u
\end{cases}
\end{align}
and proposed a time-asymptotic expansion 
\begin{align}\label{abE:2.2nm}
\begin{cases}
&\tilde v_k=\bar v+\sum\limits_{i=1}^{k}(1+t)^{-i\sigma}v_i(\xi),\\
&\tilde u_k=\bar u+\sum\limits_{i=1}^{k} (1+t)^{-(i+\frac{1}{2})\sigma}u_i(\xi), \quad \xi=\frac{x}{(1+t)^{\frac{1+\lambda}{2}}}, \,\sigma=1-\lambda.
\end{cases}
\end{align}
If the expansion holds, then the asymptotic behavior of solutions beyond the diffusion wave $\bar v(\xi)$ can be  accurately characterized.
We justified the expansion as $\lambda \in (\frac17,1)$ and further conjectured the expansion still holds for any $\lambda \in [0,1)$ in \cite{Geng-Huang-Jin-Wu}.

 The aim of this paper is to justify the expansion \eqref{abE:2.2nm} as $\lambda=0$, i.e., the constant damping system \eqref{bE:2.1}. We consider the case of $k=1$, namely,
\begin{align}\label{bE:2.2}
\begin{cases}
&v_{*}(x,t)=\bar v+(1+t)^{-1}v_1(\xi),\\
&u_{*}(x,t)=\bar u+ (1+t)^{-\frac{3}{2}}u_1(\xi),
\end{cases}
\end{align}
where $\xi=\frac{x}{\sqrt{1+t}}$, $(\bar v,\bar u)$ is the diffusion wave of \eqref{E:1.3}-\eqref{E:2.1new} and the subsequent term $(v_1, u_1)(\xi)$ 
 will be determined in section \ref{section2} below.
 
Without loss of generality, we focus on the case of $u_+=u_-=0$. The other cases can be treated by introducing a correction function 
\begin{align}
\begin{cases}
\hat v(x,t)=-(u_+-u_-)m_0(x)e^{-t},\\
\hat u(x,t)=e^{-t}
\left(u_-+(u_+-u_-)\int_{-\infty}^xm_0(y)dy\right),
\end{cases}
\end{align}
where $m_0(x)$ is a smooth function with compact support satisfying $\int_{\mathbb {R}}m_0(x)dx=1$, see \cite{Hsiao-Liu-1} for the details.
Let
\begin{equation}\label{p.n}
V(x,t)=\int_{-\infty}^{x}\left(v(y,t)-v_{*}\left(\frac{y}{\sqrt{1+t}},t\right)\right)dy,\quad V_1(x,t)= u(x,t)- u_{*}(\xi,t),
\end{equation}
and $V_0(x)=:V(x,0), V_1(x)=:V_1(x,0)$.
Assume that the initial data 
\begin{equation}\label{12}
(V_0,V_1)(x)\in H^{5}(\mathbb{R})\times H^{4}(\mathbb{R}), \quad \tilde V_0(x),\tilde V_1(x)\in L^{2}(\mathbb{R})\cap L^{1}(\mathbb{R})
\end{equation}
 with $V_0(x)=\partial_{x}\tilde V_0(x),\, V_1(x)=\partial_{x}\tilde V_1(x)$.
Set
\[N(0)=:\|V_0\|_{5}+\|V_1\|_{4}, \quad \delta=:|v_+-v_-|, \quad \delta_1=:\|\tilde V_0\|_{L^1}+\|\tilde V_1\|_{L^1}.\]
Then we have
\begin{theorem}\label{theorem 1.1}
Suppose that the initial data $V_0(x),V_1(x)$ satisfies \eqref{12}.
Then there is a small constant $\epsilon_0>0$ such that if 
\[\sqrt{N^2(0)+\delta}+\delta_1 \leq\epsilon_0,\]
the Cauchy problem \eqref{bE:2.1}-\eqref{E:1.2} admits a unique smooth solution $(v, u)(x,t)$. Moreover, it holds that
\begin{align}
\sum_{l+k\leq 5, l\leq 1}(1+t)^{l+\frac{k}{2}+\frac{3}{4}}\|\partial_t^l\partial_x^k V(\cdot, t)\|_{L^2}\leq C(\sqrt{N^2(0)+\delta}+\delta_1).
\end{align}
\end{theorem}
Furthermore, noting the relationship between $v-v_*, u-u_*$ and $V$ in Section \ref{section3} below,  we use the Cauchy-Schwartz inequality to get the following Corollary.
\begin{corollary}\label{corollary1.1}
Under the assumptions of Theorem \ref{theorem 1.1}, it holds that 
\begin{align}\label{3E:3.19}
\|(v-v_*, u-u_*)(t)\|_{L^{\infty}}\leq C(\sqrt{N^2(0)+\delta}+\delta_1) ((1+t)^{-\frac{3}{2}}, (1+t)^{-2}),
\end{align}
which justifies the expansion \eqref{bE:2.2}.
\end{corollary}

We now sketch the main strategy. 
 To justify the time asymptotic expansion \eqref{bE:2.2},
it remains to show that the remainder $v-v_*=:V_x$ decays faster than $(1+t)^{-1}$. To this end,
we firstly reduce the perturbation system into a nonlinear wave equation 
\begin{align}\label{o1}
V_{tt}+(P^{\prime}(v_{*})V_x)_x+V_t=g_{1x}+S[v_*]
\end{align}
and establish the basic estimates for $V_{tt}$ by energy method, where $g_1$ is the nonlinear term and $S[v_*]$ is the error term induced by the expansion \eqref{bE:2.2}, see \eqref{2.6} below. Motivated by the diffusion phenomenon observed in \cite{Hsiao-Liu-1} and \cite{Nishihara-Wang-Yang}, we then
regard \eqref{o1} as a diffusion equation with source term, i.e.,
\begin{align}\label{oq1}
 V_t+(P^{\prime}(v_{*})V_x)_x=g_{1x}+S[v_*]-V_{tt},
\end{align}
and use the approximate Green function $G(x,t;y,s)$ to obtain the integral formula of $V(x,t)$ through Duhamel principle,
\begin{align}\label{bvE:3.5}
V(x,t)=&\int_{R} G(x,t;y,0)V_0(y)dy+\int_0^t\int_{\mathbb R} G(x,t;y,s)[g_{1y}+S[v_*]-V_{ss}]dyds
\notag\\
&+\int_0^t\int_{\mathbb R} R_{G}(x,t;y,s)V(y,s)dyds, 
\end{align}
where $R_G$ denotes the difference between the Green function and approximate one $G(x,t;y,s)$.
Thanks to the expansion \eqref{bE:2.2}, $S[v_*]$ decays faster than that of \cite{Nishihara-Wang-Yang}. Thus we can obtain a better decay rate for 
$$
\int_0^t\int_{\mathbb R} G(x,t;y,s)S[v_*]dyds.
$$
On the other hand, we have new estimates for 
$$\int_0^t\int_{\mathbb R} G(x,t;y,s)g_{1y}dyds\quad \mbox{and}\quad \int_0^t\int_{\mathbb R} R_{G}(x,t;y,s)V(y,s)dyds.$$
Finally we can obtain the desired decay rate  $(1+t)^{-\frac32}$ for $V_x$, as shown in  \eqref{3E:3.19}, with the help of the energy estimates for the wave equation \eqref{o1}.
%
%
%

The arrangement of the present paper is as follows. In section \ref{section2}, we introduce the time asymptotic expansion $(v_*,u_*)(x,t)$ for \eqref{bE:2.1}-\eqref{E:1.2}. In section \ref{section3}, we reduce the perturbation system to a nonlinear wave equation and establish the basic estimates by energy method. The section \ref{section4} is devoted to proving Theorem \ref{theorem 1.1} by the approximate Green function method with the help of the basic energy estimates.

\

\noindent{\bf Notations}. \ \ Throughout this paper the symbol $c,\,C$
 will be used to represent a generic constant which is
independent of $x$ and $t$ and may vary from line to line.
$\|\cdot\|_{L^p}$ stands for the
$L^p(\mathbb{R})$-norm $(1\leq p\leq \infty)$. The $L^2$-norm on $\mathbb{R}$ is simply
denoted by $\|\cdot\|$. Moreover, the domain $\mathbb{R}$ will be
often abbreviated without confusions.

\section{The time asymptotic expansion}\label{section2}
We first list some properties on the diffusion wave $\bar v(\xi)$ with $\xi=\frac{x}{\sqrt{1+t}}$ of PME \eqref{E:1.3}-\eqref{E:2.1new}.

 \begin{lemma}[\cite{Duyn-Peletier-1977, Hsiao-Liu-2}]\label{lamm-01}
For the diffusion wave ${\bar v}(\xi)$  of PME \eqref{E:1.3}-\eqref{E:2.1new}, it holds that
 \begin{align}\label{bE:2.005}
&|{\bar v}(\xi)-v_+|_{\xi>0}+|\bar v(\xi)-v_-|_{\xi<0}\leq  C(|v_+-v_-|) e^{-c\xi^2},\\
&|\partial_{x}^k\partial_{t}^l\bar v|\leq C(|v_+-v_-|)(1+t)^{-\frac{k}{2}-l}e^{-c\xi^2},\hspace{0.2cm}k+l\geq 1, k,l\geq 0, \label{bE:2.4}\\
&\|\partial_{x}^k\partial_{t}^l\bar v\|\leq C(|v_+-v_-|)(1+t)^{-\frac{k}{2}-l+\frac14},\hspace{0.2cm}k+l\geq 1.
\end{align}
 \end{lemma}

As in \cite{Geng-Huang-Jin-Wu}, we introduce the time asymptotic expansion 
\begin{align}\label{bE:2.2n}
\begin{cases}
&v_{*}(x,t)=\bar v+(1+t)^{-1}v_1(\xi),\\
&u_{*}(x,t)=\bar u+ (1+t)^{-\frac{3}{2}}u_1(\xi),
\end{cases}
\end{align}
where $\bar u$ is given in \eqref{a2} and $(v_1,u_1)(\xi)$ will be determined as follows. 
Note that $\bar v_t-\bar u_x=0$, we expect 
\begin{equation}\label{2.4}
((1+t)^{-1}v_1(\xi))_t-((1+t)^{-\frac{3}{2}}u_1(\xi))_x=0,\quad \mbox{equivalently,}\quad v_{*t}-u_{*x}=0,
\end{equation}
which implies that
\begin{equation}\label{vE:2.5}
u_1(\xi)=-\frac{1}{2}\xi v_1-\frac{1}{2}\int_{-\infty}^{\xi}v_1d\xi=-\frac12(\xi G_1)_{\xi},
\end{equation}
where $G_1=\int_{-\infty}^{\xi}v_1d\xi$. Denote the source term
\begin{align}\label{2.6}
S[v_*]=u_{*t}+P(v_{*})_x+u_{*}.
\end{align}
For convenience, we also write $P^{\prime}(v)=\frac{dP(v)}{dv}$. From Darcy's law $P(\bar v)_{x}+\bar{u}=0$ and \eqref{vE:2.5}, a straightforward computation shows that
\begin{align}\label{q}
S[v_*]&=\bar u_t+((1+t)^{-1}P^{\prime}(\bar v)v_1)_x+(1+t)^{-\frac{3}{2}}u_1+g_{2x}+g_{3t}
\notag\\
&=(1+t)^{-\frac32}\Big(P^{\prime}(\bar v)G_{1\xi}-\frac{1}{2}\xi G_1+\frac{1}{2}\xi P(\bar v)_{\xi}\Big)_{\xi}+g_{2x}+g_{3t}
\notag\\
&=(1+t)^{-\frac32}\Big(P^{\prime}(\bar v)G_{1\xi}-\frac{1}{2}\xi G_1+\frac{1}{2}\xi P(\bar v)_{\xi}\Big)_{\xi}+O(1)(1+t)^{-\frac52},
\end{align}
where 
\begin{align}
g_2&=P(v_{*})-P(\bar v)-P^{\prime}(\bar v)(v_{*}-\bar v),\label{vE:2.9}\\
g_3&=(1+t)^{-\frac{3}{2}}u_1.\label{vE:2.10}
\end{align}
To eliminate the term concerning $(1+t)^{-\frac32}$ in \eqref{q}, we seek for $G_1(\xi)$ satisfying
\begin{align}\label{ode}
P^{\prime}(\bar v)G_{1\xi}-\frac{1}{2}\xi G_1+\frac{1}{2}\xi P(\bar v)_{\xi}=0
\end{align}
and $v_1(\xi)=G_{1{\xi}}$,
Thus, $
S[v_*]=g_{2x}+g_{3t}.
$
We choose a solution $G_1$ to \eqref{ode} by
\begin{equation}\label{E:2.13}
G_1(\xi)=-e^{\int \frac{\xi}{2P^{\prime}(\bar v)}d\xi} \int_0^{\xi}e^{-\int \frac{\eta}{2P^{\prime}(\bar v)}d\eta}\frac{\eta \bar v^{\prime}_{\eta}}{2}d\eta.
\end{equation}
\begin{lemma}
Let $\delta=|v_+-v_-|$. Then for $(v_1, u_1)$ given in \eqref{vE:2.5} and \eqref{E:2.13}, it holds that 
\begin{align}
S[v_*]=O(1)\delta e^{-c\xi^2}(1+t)^{-\frac52}.
\end{align}
\end{lemma}
\begin{proof}
It follows from \eqref{vE:2.9}-\eqref{vE:2.10} that 
\[g_2=O(1)(1+t)^{-2}v_1^2,\quad g_3=(1+t)^{-\frac{3}{2}}u_1.\]
By Lemma \ref{lamm-01}, \eqref{vE:2.5} and \eqref{E:2.13}, we have
\[S[v_*]=g_{2x}+g_{3t}=O(1)\delta e^{-c\xi^2}(1+t)^{-\frac52}.\]
Thus the proof is completed.
\end{proof}

Moreover, we have 
\begin{lemma}\label{blemma2.2}
It holds that 
\begin{align}
&\|g_2\|_{L^1}\leq C\delta (1+t)^{-\frac{3}{2}}, \quad  \|\partial_t^n\partial_x^k g_2\|_{L^2}\leq C\delta (1+t)^{-n-\frac{k}{2}-\frac{7}{4}},\label{1vE:2.16}\\
&\|g_3\|_{L^1}\leq C\delta (1+t)^{-1}, \quad  \|\partial_t^n\partial_x^k g_3\|_{L^2}\leq C\delta (1+t)^{-n-\frac{k}{2}-\frac{5}{4}}.\label{1vE:2.17}
\end{align}\end{lemma}
\begin{proof}
Lemma \ref{blemma2.2} can be proved by the straightforward computations on \eqref{vE:2.9}-\eqref{vE:2.10}.
\end{proof}

\section{reduced system}\label{section3}
From \eqref{bE:2.1}, \eqref{2.4} and \eqref{2.6} we get the following perturbation system:
\begin{align}\label{vE:2.13}
\begin{cases}
(v-v_{*})_t-(u-u_{*})_x=0,\\
(u-u_{*})_t+(P^{\prime}(v_{*})(v-v_{*}))_x+u-u_{*}=g_{1x}+S[v_*]
\end{cases}
\end{align}
with the initial data 
\begin{align}\label{3.2}
(v-v_{*},u-u_{*})(x,0)=(v_0(x)-v_*(x,0),u_0(x)-u_*(x,0)),
\end{align}
where  
\begin{align}\label{vE:2.8}
g_1=-(P(v)-P(v_{*})-P^{\prime}(v_{*})( v-v_{*})).
\end{align}
Let
\begin{equation*}
V=\int_{-\infty}^x(v(y,t)-{v_{*}}(y,t)) dy,
\end{equation*}
which yields that $V_x=v-v_*$ and $V_t=u-u_*$.
Then we can rewrite \eqref{vE:2.13}-\eqref{3.2} as a nonlinear wave equation 
\begin{align}\label{1E:2.19}
V_{tt}+(P^{\prime}(v_{*})V_x)_x+V_t=g_{1x}+S[v_*]
\end{align}
with the initial data
\begin{align}\label{1E:2.191}
(V,V_t)(x,0)=(V_0, V_1)(x).
\end{align}

\
Following the framework of \cite{Nishihara}, we seek for the solution $V(x,t)$ of \eqref{1E:2.19} in the following solution space
\[X_T=:\{V(x,t)\,|\,V \in C([0,T);H^{5}(\mathbb R)), V_t \in C([0,T);H^{4}(\mathbb R))\}.\]
Since the local existence of the
solution of \eqref{1E:2.19} can be proved by the standard iteration method, see \cite{Matsumura-1997}, the main effort in this section is to establish
the a priori estimates for the solution.

For any $T\in (0,+\infty)$, define
\begin{align}\label{vE:3.37}
N(T)=\sup \limits_{0\leq t\leq T}\sum\limits_{l+k\leq 5, l\leq 1}(1+t)^{l+\frac{k}{2}+\frac{3}{4}}\|\partial_t^l\partial_x^k V(t)\|_{L^2}.
\end{align}
Suppose that $N(T)\leq \epsilon$, where $\epsilon$ is sufficiently small and will be determined later. 

\begin{lemma}\label{newlemma2.1}
For any $T>0,$ assume that $V(x,t) \in X_T$ is the solution of \eqref{1E:2.19}. If $\epsilon$  and $\delta=\left|v_{+}-v_{-}\right|$ are small, then it holds that for any $0<t<T$,
 \begin{align}\label{2.20}
&\left\|V(\cdot, t)\right\|^{2}+\sum_{k=0}^{4}(1+t)^{k+1}(\left\|\partial_{x}^{k} V_t(\cdot, t)\right\|^{2} +\left\|\partial_{x}^{k} V_x(\cdot, t)\right\|^{2})
 \leq C\left(N^2(0)+\delta\right)
\end{align}
and 
 \begin{align}\label{2.21}
&(1+t)^{2}\left\|V_{t}(\cdot, t)\right\|^{2}+\sum_{k=0}^{3}(1+t)^{k+3}(\left\|\partial_{x}^{k} V_{tt}(\cdot, t)\right\|^{2} +\left\|\partial_{x}^{k} V_{xt}(\cdot, t)\right\|^{2})
 \leq C\left(N^2(0)+\delta\right).
\end{align}
\end{lemma}
\begin{proof}
The proof can be completed by the same line as in \cite{Nishihara} and the details are omitted. 
\end{proof}

The decay rates of $\partial_x^k\partial_t^lV_{tt}$ can be improved by the following estimates.
\begin{lemma}\label{newlemma2.2}
Under the assumptions of Lemma \ref{newlemma2.1}, it holds that  for any $0<t<T$,
 \begin{align}\label{2.22}
&(1+t)^{4}\left\|V_{tt}(\cdot, t)\right\|^{2}+\sum_{k=0}^{2}(1+t)^{k+5}(\left\|\partial_{x}^{k} V_{ttt}(\cdot, t)\right\|^{2} +\left\|\partial_{x}^{k} V_{xtt}(\cdot, t)\right\|^{2})
 \leq C\left(N^2(0)+\delta\right).
\end{align}
\end{lemma}
\begin{proof}
For $k=0,1,2$, taking the procedure as 
\[\int_0^t\int_{\mathbb R}\Big[\partial_x^k\partial_t^2 \eqref{1E:2.19}\times (\mu+t)^{k+4} \partial_x^k\partial_t^3V+\partial_x^k\partial_t^2 \eqref{1E:2.19}\times (\mu+t)^{k+4} \partial_x^k\partial_t^2V\Big]dxd\tau\]
 and choosing constant $\mu$ large enough, together with Lemmas \ref{lamm-01} and \ref{newlemma2.1}, yield that
  \begin{align}\label{2.23}
 &(1+t)^{4}\left\|V_{tt}(\cdot, t)\right\|^{2} +(1+t)^{k+4}(\left\| \partial_x^kV_{ttt}(\cdot, t)\right\|^{2} +\left\|\partial_x^kV_{xtt}(\cdot, t)\right\|^{2})
\notag\\
+&\int_{0}^{t}(1+\tau)^{k+4}\left[\left\|\partial_x^kV_{xtt}(\cdot, \tau)\right\|^{2}+\left\| \partial_x^kV_{ttt}(\cdot, \tau)\right\|^{2}\right] d \tau  \leq C\left(N^2(0)+\delta\right).
  \end{align}
In the similar way as for
\[\int_0^t\int_{\mathbb R}\partial_x^k\partial_t^2 \eqref{1E:2.19}\times (\mu+t)^{k+5} \partial_x^k\partial_t^3Vdxd\tau,\]
we have
  \begin{align}\label{2.24}
  (1+t)^{k+5}(\left\| \partial_x^kV_{ttt}(\cdot, t)\right\|^{2} +\left\|\partial_x^kV_{xtt}(\cdot, t)\right\|^{2})
+\int_{0}^{t}(1+\tau)^{k+5}\left\| \partial_x^kV_{ttt}(\cdot, \tau)\right\|^{2} d \tau 
 \leq C\left(N^2(0)+\delta\right).
  \end{align}
 Thus, \eqref{2.23} and \eqref{2.24} lead to \eqref{2.22}. Therefore, the proof of Lemma \ref{newlemma2.2} is completed.
 \end{proof}
 Similarly, we have
\begin{lemma}\label{newlemma2.3}
Under the assumptions of Lemma \ref{newlemma2.1}, it holds that  for any $0<t<T$,
 \begin{align}\label{2.25}
&(1+t)^{6}\left\|V_{ttt}(\cdot, t)\right\|^{2}+\sum_{k=0}^{1}(1+t)^{k+7}(\left\|\partial_{x}^{k} V_{tttt}(\cdot, t)\right\|^{2} +\left\|\partial_{x}^{k} V_{xttt}(\cdot, t)\right\|^{2})
 \leq C\left(N^2(0)+\delta\right).
\end{align}
\end{lemma}
\begin{lemma}\label{newlemma2.4}
Under the assumptions of Lemma \ref{newlemma2.1}, it holds that  for any $0<t<T$,
 \begin{align}\label{2.26}
&(1+t)^{8}\left\|\partial_t^4V(\cdot, t)\right\|^{2}+(1+t)^{9}(\left\|\partial_{t}^{5} V(\cdot, t)\right\|^{2} +\left\|\partial_{t}^{4} V_{x}(\cdot, t)\right\|^{2})
 \leq C\left(N^2(0)+\delta\right).
\end{align}
\end{lemma}

\section{Green function method}\label{section4}
Note that the decay rates obtained in Lemmas \ref{newlemma2.1}-\ref{newlemma2.4} are not fast enough to close the a priori assumption $N(T)\leq \epsilon$ in \eqref{vE:3.37}.
In this section, we will use the approximate Green function method to improve the decay rates so that the a priori assumption can be closed. Since Lemmas \ref{newlemma2.1}-\ref{newlemma2.4} have provided the desired estimates to close the a priori assumption for local time, we focus on the large time $t>1$ in what follows.
As in \cite{Nishihara-Wang-Yang}, we rewrite \eqref{1E:2.19} as
\begin{align}\label{1E:2.21}
V_t+(a(x,t)V_x)_x=g_{1x}+S[v_*]-V_{tt},
\end{align}
where $a(x,t)=P^{\prime}(v_{*})$ and construct a minimizing Green function as 
\begin{align*}
G(x,t;y,s)=\Big(\frac{-1}{4\pi a(x,t)(t-s)}\Big)^{\frac{1}{2}}\exp \Big(\frac{(x-y)^2}{4A(y,s;t)(t-s)}\Big)
\end{align*}
satisfying the basic requirement 
\begin{equation}\label{vE:3.3}
G(x,t;y,t)=\delta (y-x),
\end{equation}
where $\delta$ is the Dirac function and $A(y,s;t)=P^{\prime}(v_{*}(\eta,t))$ and 
\begin{align*}
\eta=
\begin{cases}
y/\sqrt{1+s}, \quad s>t/2,\\
y/\sqrt{1+t/2}, \quad s\leq t/2.
\end{cases}
\end{align*}
Recall that $S[v_*]=g_{2x}+g_{3t}$, then the solution $V(x,t)$ to \eqref{1E:2.21} can be written as the integral form
\begin{align}\label{vE:3.5}
V(x,t)=&\int_{R} G(x,t;y,0)V_0(y)dy+\int_0^t\int_{\mathbb R}G(x,t;y,s)[g_{1y}+g_{2y}+g_{3s}-V_{ss}]dyds
\notag\\
&+\int_0^t\int_{\mathbb R}R_{G}(x,t;y,s)V(y,s)dyds,
\end{align}
where 
\begin{equation*}
 R_{G}(x,t;y,s)=G_s(x,t;y,s)-\Big\{a(y,s)G_y(x,t;y,s)\Big\}_y.
\end{equation*}

\subsection{Properties of approximate Green function}
In this subsection, we recall the properties of $G(x,t;y,s)$ and $R_{G}(x,t;y,s)$ introduced in \cite{Nishihara-Wang-Yang}.
\begin{lemma}[\cite{Nishihara-Wang-Yang}]\label{vlemma3.1}
For $l, h\leq 1$, it holds that 
\begin{align}\label{vE:3.8}
|\partial_t^l\partial_s^h\partial_x^k\partial_y^mG(x,t;y,s)|=&O(1)\Big(\sum\limits_{m_1+m_2=m}(t-s)^{-\frac{m_1}{2}}\theta^{m_2}\Big)\Big(\sum\limits_{k_1+k_2=k}(t-s)^{-\frac{k_1}{2}}(1+t)^{-\frac{k_2}{2}}\Big)
\notag\\
&\times ((1+t)^{-1}+(t-s)^{-1})^l(\theta_1^2+(t-s)^{-1})^h G_D(x-y, t-s),
\end{align}
 where 
\begin{equation*}
G_D(y,s)=\Big(\frac{1}{4\pi s}\Big)^{\frac{1}{2}}exp\Big(\frac{-y^2}{Ds}\Big)\quad \mbox{with}\quad \| G_D\|_{L^{p}}\leq Cs^{-\frac{1}{2}(1-\frac{1}{p})}\quad \mbox{for}\quad p \geq 1,
\end{equation*}
and $\theta(s)=\theta_1(s)+\theta_2(s)$ with 
\begin{align*}
\theta_1(s)=
\begin{cases}
(1+s)^{-\frac{1}{2}}, \quad s>t/2,\\
0, \quad s\leq t/2,
\end{cases}
\quad \mbox{and }\quad 
\theta_2(s)=
\begin{cases}
0, \quad s>t/2,\\
(1+t)^{-\frac{1}{2}}, \quad s\leq t/2.
\end{cases}
\end{align*}
\end{lemma}

\begin{lemma}[\cite{Nishihara-Wang-Yang}]\label{vlemma3.2}
It holds that
\begin{equation}\label{vE:3.9}
R_G=O(1)\delta \Theta(t,s)\tilde E(y,t,s)G_D(x-y, t-s), 
\end{equation}
where 
\begin{align*}
\Theta(t,s)&=
\begin{cases}
(1+s)^{-1}+(t-s)^{-\frac{1}{2}}(1+s)^{-\frac{1}{2}},\quad s>t/2, \\
(1+t)^{-1}+(t-s)^{-\frac{1}{2}}(1+s)^{-\frac{1}{2}}, \quad s\leq t/2,
\end{cases}\\
\tilde E(y,t,s)&=
\begin{cases}
E(y,s), \quad s>t/2,\\
E(y,t), \quad s\leq t/2,
\end{cases}
\end{align*}
with
$E(y,\tau)=exp(-Cy^2/(1+\tau))$ for some constant $C>0$. Moreover, it holds for $s<t/2$ that
\begin{align}\label{vE:3.10}
|\partial_t^l\partial_x^kR_G(x,t;y,s)|\leq C\delta (1+s)^{-\frac{1}{2}}(t-s)^{-l-\frac{k+1}{2}}E(y,t)G_D(x-y, t-s),
\end{align}
and for $s=t/2$ that
\begin{align}\label{vE:3.11}
\lim \limits_{s\rightarrow t/2^{\pm}}|\partial_x^kR_G(x,t;y,s)|\leq C\delta t^{-1-\frac{k}{2}}E(y,t/2)G_D(x-y, t/2).
\end{align}
\end{lemma}
\begin{proof}
Although Lemma \ref{vlemma3.2} was obtained in \cite{Nishihara-Wang-Yang}, we give more details to obtain \eqref{vE:3.9} and some details will be used later.
The direct computation shows that 
\begin{align}\label{vE:3.15}
R_G=&\frac{G}{2(t-s)}\Big(1-\frac{a(y,s)}{A(y,s;t)}\Big)+\frac{(x-y)^2G}{4A(y,s)(t-s)^2}\Big(1-\frac{a(y,s)}{A(y,s;t)}\Big)+G \tilde R,
\end{align}
where $\tilde R$ satisfies 
\begin{align}
\tilde R&=O(1)\delta
\begin{cases}
((1+s)^{-1}+(t-s)^{-\frac{1}{2}}(1+s)^{-\frac{1}{2}}) E(y,s), \quad s>t/2,\\
(t-s)^{-\frac{1}{2}}(1+s)^{-\frac{1}{2}}E(y,s)+(t-s)^{-\frac{1}{2}}(1+t)^{-\frac{1}{2}}E(y,t),\quad s\leq t/2\label{vvE:3.16}
\end{cases}
\\
&=O(1)\delta \Theta(t,s) E(y,t).\nonumber
\end{align}
When $s>t/2$, it follows from $a(y,s)=A(y,s;t)$ that $R_G=G\tilde R$, which 
 implies that \eqref{vE:3.9} holds.
It remains to show \eqref{vE:3.9} for $s\leq \frac{t}{2}$. The straightforward computations show that 
\begin{align}\label{vE:3.16}
&\frac{1}{(t-s)}\Big|1-\frac{a(y,s)}{A(y,s;t)}\Big|=\Big|\frac{1}{t-s}\frac{p^{\prime}(v_{*}(\frac{y}{\sqrt{1+t/2}}))-p^{\prime}(v_{*}(\frac{y}{\sqrt{1+s}}))}{p^{\prime}(v_{*}(\frac{y}{\sqrt{1+t/2}}))}\Big|
\notag\\
=&\Big|\frac{1}{t-s}\frac{p^{\prime\prime}(v_{*}(\frac{y}{\sqrt{1+\theta}}))v_{*}^{\prime}(\frac{y}{\sqrt{1+\theta}})}{p^{\prime}(v_{*}(\frac{y}{\sqrt{1+t/2}}))}\Big(\frac{y}{\sqrt{1+t/2}}-\frac{y}{\sqrt{1+s}}\Big)\Big|
\leq
O(1)\delta E(y,t)(t-s)^{-\frac{1}{2}}(1+s)^{-\frac{1}{2}},
\end{align}
which, together with \eqref{vvE:3.16}, leads to \eqref{vE:3.9} directly.
\end{proof}

%
%
%
Next, we will summarize the properties of $G(x,t;y,s)$ and $R_{G}(x,t;y,s)$ introduced in \cite{Nishihara-Wang-Yang}.
\begin{lemma}\label{vlemma3.3}
For $s>\frac{t}{2}$ and $k_1\geq 1, k_2\geq 0$, it holds that 
\begin{align}
&\partial_x^{k_1}G=(-1)^{k_1}\partial_y^{k_1}G+\sum\limits_{0\leq \beta<{k_1}}C_{\beta}\partial_y^{\beta}(G\tilde R^1_{{k_1}-\beta-1}), \label{vE:3.20}\\
&\partial_t\partial_x^{k_2}G=-\partial_s\partial_x^{k_2}G+\sum\limits_{0\leq \beta\leq k_2}C_{\beta}\partial_y^{\beta}(G\tilde h^1_{k_2-\beta}),\label{vE:3.21}
\end{align}
where $\tilde R^1_l$ and $\tilde h^1_l$ represent  some  generic functions satisfying 
\begin{equation}\label{vE:3.23}
\tilde R^1_l=O(1)\delta (1+t)^{-\frac{l+1}{2}}E(y,t)\quad \mbox{and}\quad  \partial_y\tilde R^1_l=O(1)\delta (1+t)^{-\frac{1+l}{2}}[(t-s)^{-\frac{1}{2}}+(1+t)^{-\frac{1}{2}}]
\end{equation}
and 
\begin{equation}\label{1vE:3.23}
\tilde h^1_l=O(1)\delta (1+t)^{-1-\frac{l}{2}}E(y,t)\quad \mbox{and}\quad \partial_y\tilde h^1_l=O(1)\delta (1+t)^{-1-\frac{l}{2}}[(t-s)^{-\frac{1}{2}}+(1+t)^{-\frac{1}{2}}], \quad \forall\, l\geq 0.
\end{equation}
Furthermore, for the above mentioned $\tilde R^1_l$, it holds that
\begin{equation}\label{11vE:3.23}
\partial_x(G\tilde R^1_l)=-\partial_y(G\tilde R^1_{l})+G\tilde R^1_{l+1}.
\end{equation}
\end{lemma}
\begin{proof}
Lemma \ref{vlemma3.3} can be proved by the straightforward computations and Lemma \ref{vlemma3.1}. The details are omitted.
\end{proof}

\begin{lemma}\label{vlemma3.4}
For $s>\frac{t}{2}$ and $k_1\geq 1, k_2\geq 0$, it holds that 
\begin{align}
&\partial_x^{k_1}R_G=(-1)^{k_1}\partial_y^{k_1}R_G+\sum\limits_{\beta<{k_1}}C_{\beta}\partial_y^{\beta}(G\tilde R^2_{{k_1}-\beta-1}), \label{vE:3.28}\\
&\partial_t\partial_x^{k_2}R_G=-\partial_s\partial_x^{k_2}R_G+\sum\limits_{\beta\leq {k_2}}C_{\beta}\partial_y^{\beta}(G\tilde h^2_{{k_2}-\beta}), \label{vE:3.29}
\end{align}
where $\tilde R^2_l$ and $\tilde h^2_l$  represent  some  generic functions satisfying 
\begin{equation}\label{vE:3.31}
\tilde R^2_l=O(1)\delta \Theta(t,s) E(y,t) (1+t)^{-\frac{1+l}{2}} \quad \mbox{and}\quad \tilde h^2_l=O(1)\delta \Theta(t,s) E(y,t)(1+t)^{-1-\frac{l}{2}},\quad \forall\, l\geq 0.
\end{equation}
\end{lemma}
\begin{proof}
Recall from \eqref{vvE:3.16} that for $ s>\frac{t}{2}$, $R_G=G\tilde R$ and $\tilde R=O(1)\delta \Theta(t,s) E(y,t)$ satisfies
\begin{align*}
\partial_x \tilde R&=-\partial_y \tilde R+\hat R \quad \mbox{with}\quad \hat R=O(1)\delta \Theta(t,s) E(y,t)(1+t)^{-\frac{1}{2}},\notag\\
\partial_t \tilde R&=-\partial_s \tilde R+\hat h\quad \mbox{with}\quad \hat h=O(1)\delta \Theta(t,s)E(y,t) (1+t)^{-1}.
\end{align*}
Then by the direct computations and Lemma \ref{vlemma3.3}, we can verify 
\eqref{vE:3.28}-\eqref{vE:3.29}. Thus, the proof of Lemma \ref{vlemma3.4} is completed.
\end{proof}

Note that \eqref{vE:3.10} gives
\begin{align*}
\|\partial_t^l\partial_x^k R_G\|_{L^1}
= O(1)\delta t^{-l-\frac{k+1}{2}}(1+s)^{-\frac{1}{2}} \quad \mbox{for}\quad s\leq \frac{t}{2},
\end{align*}
whose decay rate is not fast enough in the analysis. We need to derive a better decay rate as follows.
\begin{lemma}\label{vcorollary3.1}
For $s< t/2$, it holds that
\begin{align*}
\|\partial_t^l\partial_x^k R_G\|_{L^1}= & O(1)\delta [t^{-l-\frac{k+1}{2}-\frac{1}{2q_1}}(1+s)^{-\frac{1}{2}+\frac{1}{2q_1}}
+ t^{-l-\frac{k}{2}-1}]\quad \mbox{for}\quad q_1\geq 1.
\end{align*}
\end{lemma}
\begin{proof}
Different from \eqref{vE:3.16}, we estimate $1-\frac{a(y,s)}{A(y,s;t)}$ by
\begin{align*}
&1-\frac{a(y,s)}{A(y,s;t)}=\frac{p^{\prime}(v_{*}(\frac{y}{\sqrt{1+t/2}}))-p^{\prime}(v_{*}(\frac{y}{\sqrt{1+s}}))}{p^{\prime}(v_{*}(\frac{y}{\sqrt{1+t/2}}))}
\notag\\
=&O(1)\delta \left(\frac{t}{2}-s\right)\int_0^1 E(y,s+\theta(t/2-s))(1+s+\theta(t/2-s))^{-1}d\theta=:O(1)\delta \left(\frac{t}{2}-s\right) \hat F.
\end{align*}
Differentiating \eqref{vE:3.15}, using Lemma \ref{vlemma3.1} and \eqref{vvE:3.16} gives that  
\begin{align*}
\partial_t^l\partial_x^k R_G=&O(1)\delta (t-s)^{-l-\frac{k+1}{2}}\left[(1+s)^{-\frac{1}{2}}E(y,s)+(1+t)^{-\frac{1}{2}}E(y,t)+(t-s)^{-\frac12}\left(\frac{t}{2}-s\right)\hat F\right] G_D,
\end{align*}
and then
\begin{align*}
\|\partial_t^l\partial_x^k R_G\|_{L^1}
= & O(1)\delta (t-s)^{-l-\frac{k+1}{2}}\Big((1+s)^{-\frac{1}{2}}\|E(y,s)\|_{L^{q_1}}\|G_D\|_{L^{p_1}}+(1+t)^{-\frac{1}{2}}\|E(y,t)\|_{L^{2}}\|G_D\|_{L^2}\Big)
\notag\\
&+O(1)\delta (t-s)^{-l-\frac{k}{2}-1}(t/2-s)\|\hat F\|_{L^{q_2}}\| G_D\|_{L^{p_2}},
\end{align*}
where $\frac{1}{p_i}+\frac{1}{q_i}=1, p_i,q_i\geq 1$.
Thanks to the Minkowski inequality, we have
\begin{align*}
\|\hat F\|_{L^{q}}&=\Big(\int_{\mathbb R}\Big(\int_0^1 E(y,s+\theta(t/2-s))(1+s+\theta(t/2-s))^{-1}d\theta\Big)^qdy\Big)^{\frac{1}{q}}
\notag\\
&\leq \int_0^1 \Big(\int_{\mathbb R}[E(y,s+\theta(t/2-s))(1+s+\theta(t/2-s))^{-1} ]^qdy\Big)^{\frac{1}{q}} d\theta
\notag\\
&\leq \int_0^1 \Big(\int_{\mathbb R}e^{-\frac{Cqy^2}{1+s+\theta(t/2-s)}}dy\Big)^{\frac{1}{q}} (1+s+\theta(t/2-s))^{-1} d\theta
\notag\\
&= O(1)\int_0^1  (1+s+\theta(t/2-s))^{-1+\frac{1}{2q}} d\theta
\notag\\
&= O(1)(t/2-s)^{-1}(1+t)^{\frac{1}{2q}},
\end{align*}
which leads to 
\begin{align*}
\|\partial_t^l\partial_x^k R_G\|_{L^1}
= & O(1)\delta (t-s)^{-l-\frac{k+1}{2}-\frac{1}{2}(1-\frac{1}{p_1})}(1+s)^{-\frac{1}{2}+\frac{1}{2q_1}}+O(1)\delta (t-s)^{-l-\frac{k+1}{2}-\frac{1}{4}}(1+t)^{-\frac{1}{4}}
\notag\\
&+O(1)\delta (t-s)^{-l-\frac{k}{2}-1-\frac{1}{2}(1-\frac{1}{p_2})}(1+t)^{\frac{1}{2q_2}}
\notag\\
= & O(1)\delta t^{-l-\frac{k+1}{2}-\frac{1}{2q_1}}(1+s)^{-\frac{1}{2}+\frac{1}{2q_1}}
+O(1)\delta t^{-l-\frac{k}{2}-1}.
\end{align*}
Thus, the proof of Lemma \ref{vcorollary3.1} is completed.
\end{proof}

\subsection{Decay rate}
In this subsection, we will derive the sharper decay rate for $V(x,t)$ to close the a priori assumption $N(T)\leq \epsilon$ in \eqref{vE:3.37}. From \eqref{vE:3.5} we have
\begin{equation}\label{4.43}
\partial_t^l\partial_x^k V(x,t)=\sum_{i=1}^6I^{l,k}_i,
\end{equation}
where $l\leq 1, k+l\leq 5$ and 
\begin{align*}
I^{l,k}_1&=\int_{\mathbb R}\partial_t^l\partial_x^k G(x,t;y,0)V(y,0)dy,
\notag\\
I^{l,k}_2&=\partial_t^l\int_0^t\int_{\mathbb R}\partial_x^kG(x,t;y,s)g_{2y}dyds,
\notag\\
I^{l,k}_3&=\partial_t^l\int_0^t\int_{\mathbb R}\partial_x^kG(x,t;y,s)g_{3s}dyds,
\notag\\
I^{l,k}_4&=\partial_t^l\int_0^t\int_{\mathbb R}\partial_x^kG(x,t;y,s)g_{1y}dyds,
\notag\\
I^{l,k}_5&=-\partial_t^l\int_0^t\int_{\mathbb R}\partial_x^kG(x,t;y,s)V_{ss}dyds,
\notag\\
I^{l,k}_6&=\partial_t^l\int_0^t\int_{\mathbb R}\partial_x^kR_{G}(x,t;y,s)Vdyds.
\end{align*}

For $I^{l,k}_1$, we use Lemma \ref{vlemma3.1}, the facts $V(y,0)=\partial_y\tilde V(y,0)$ and $\tilde V(y,0) \in L^{2}\cap L^1$ with $\|\tilde V(y,0)\|_{L^1}<\delta_1$ to get
\begin{align}\label{1vE:3.44}
\|I^{l,k}_1\|_{L^2}&
=\|\int_{\mathbb R}\partial_t^l\partial_y \partial_x^k G(x,t;y,0)\tilde V(y,0)dy\|_{L^2}=\| \partial_t^l\partial_y\partial_x^k G(x,t;y,0)\ast \tilde V(y,0)\|_{L^2}
\notag\\
&=\| \partial_t^l\partial_y \partial_x^k G(x,t;y,0)\|_{L^2} \|\tilde V(y,0) \|_{L^1}
=O(1)\delta_1(1+t)^{-l-\frac{k}{2}-\frac{3}{4}},
\end{align}
where we have used the Young inequality $\|f\ast g\|_{L^r}\leq \|f\|_{L^p}\|g\|_{L^q}$ with $1+\frac{1}{r}=\frac{1}{p}+\frac{1}{q}$ for $1\leq p,q\leq \infty$.

For $I^{l,k}_2$ with $l=0, k\leq 5$, we have
\begin{align}\label{vE:3.66}
I^{0,k}_2&=\int_0^t\int_{\mathbb R}\partial_x^kG(x,t;y,s)g_{2y}dyds=\Big(\int_0^{\frac{t}{2}}+\int_{\frac{t}{2}}^t\Big)\int_{\mathbb R}\partial_x^kG(x,t;y,s)g_{2y}dyds=:I^{0,k}_{2,1}+I^{0,k}_{2,2}.
\end{align}
It follows from  \eqref{1vE:2.16} and \eqref{vE:3.8} that
\begin{align}\label{vE:3.67}
\|I^{0,k}_{2,1}\|_{L^2}&=\|\int_0^{\frac{t}{2}}\int_{\mathbb R}\partial_y\partial_x^kG(x,t;y,s)g_{2}dyds\|_{L^2}= O(1)\int_0^{\frac{t}{2}}t^{-\frac{k}{2}}(1+t)^{-\frac{1}{2}} \| G\|_{L^2} \|g_2 \|_{L^1}ds
\notag\\
&= O(1)\delta (1+t)^{-\frac{k}{2}-\frac{3}{4}}\int_0^{\frac{t}{2}}(1+s)^{-\frac{3}{2}}ds=O(1)\delta (1+t)^{-\frac{k}{2}-\frac{3}{4}}.
\end{align}
In addition, we use \eqref{1vE:2.16}, \eqref{vE:3.20}  and \eqref{vE:3.23} to get
\begin{align}\label{vE:3.68}
\|I^{0,k}_{2,2}\|_{L^2}&=\|\int_{\frac{t}{2}}^t\int_{\mathbb R}\Big[(-1)^{k}\partial_y^{k}G+\sum\limits_{\beta<{k}}C_{\beta}\partial_y^{\beta}(G\tilde R^1_{{k}-\beta-1})\Big]g_{2y}dyds\|_{L^2}
\notag\\
&=O(1)\int_{\frac{t}{2}}^t \|G_D\|_{L^1}\|\partial_y^{k+1}g_{2}\|_{L^2}ds+O(1)\delta \sum\limits_{\beta<{k}}(1+t)^{-\frac{k-\beta}{2}}\int_{\frac{t}{2}}^t \|G_D\|_{L^1}\|\partial_y^{\beta+1}g_{2}\|_{L^2}ds
\notag\\
&= O(1)\delta (1+t)^{-\frac{k}{2}-\frac{5}{4}}.
\end{align}
Substituting \eqref{vE:3.67}-\eqref{vE:3.68} into \eqref{vE:3.66} yields that for $k\leq 5$,
\begin{align}\label{new5}
\|I^{0,k}_{2}\|_{L^2}=O(1)\delta (1+t)^{-\frac{k}{2}-\frac{3}{4}}.
\end{align}
When $l=1, k\leq 4$, we deduce from \eqref{vE:3.20}-\eqref{vE:3.21} that 
\begin{align*}
I^{1,k}_2
&=\int_{\mathbb R}\partial_x^kG(x,t;y,\frac{t}{2})g_{2y}(\frac{t}{2})dy+\int_0^{\frac{t}{2}}\int_{\mathbb R}\partial_t\partial_y\partial_x^kG(x,t;y,s)g_{2}dyds+\int_{\frac{t}{2}}^t\int_{\mathbb R}G\partial_s\partial_y^{k+1}g_{2}dyds
\notag\\
&\quad+\sum\limits_{\beta<{k}}C_{\beta}\int_{\frac{t}{2}}^t\int_{\mathbb R}G\tilde R^1_{{k}-\beta-1}\partial_s\partial_y^{\beta+1}g_{2}dyds+\sum\limits_{\beta\leq k}C_{\beta}\int_{\frac{t}{2}}^t\int_{\mathbb R}G\tilde h^1_{k-\beta}\partial_y^{\beta+1}g_{2}dyds,
\end{align*}
which, together with \eqref{1vE:2.16}, \eqref{vE:3.8} and \eqref{vE:3.23}-\eqref{1vE:3.23}, yields that
\begin{align}\label{new6}
\|I^{1,k}_{2}\|_{L^2}=&O(1)t^{-\frac{k}{2}}\|G_D\|_{L^1}\|g_{2y}(\frac{t}{2})\|_{L^2}+O(1)t^{-\frac{3}{2}-\frac{k}{2}}\int_0^{\frac{t}{2}}\|G_D\|_{L^2}\|g_{2}\|_{L^1}ds
\notag\\
&+O(1)\int_{\frac{t}{2}}^t\|G_D\|_{L^1}\|\partial_s\partial_y^{k+1}g_{2}\|_{L^2}ds+O(1)\delta \sum\limits_{\beta<{k}}(1+t)^{-\frac{k-\beta}{2}}\int_{\frac{t}{2}}^t\|G_D\|_{L^1}\|\partial_s\partial_y^{\beta+1}g_{2}\|_{L^2}ds
\notag\\
&+O(1)\delta \sum\limits_{\beta\leq {k}}(1+t)^{-1-\frac{k-\beta}{2}}\int_{\frac{t}{2}}^t\|G_D\|_{L^1}\|\partial_y^{\beta+1}g_{2}\|_{L^2}ds=O(1)\delta (1+t)^{-\frac{k}{2}-\frac{7}{4}}.
\end{align}
In summary, we have from \eqref{new5} and \eqref{new6} that
\begin{align}\label{1vE:3.75}
\|I^{l,k}_{2}\|_{L^2}=O(1)\delta (1+t)^{-l-\frac{k}{2}-\frac{3}{4}}.
\end{align}

For $I^{l,k}_3$ with $l=0, k\leq 5$, we have
\begin{align}\label{vE:3.72}
I^{0,k}_3&=\int_0^t\int_{\mathbb R}\partial_x^kG(x,t;y,s)g_{3s}dyds=\Big(\int_0^{\frac{t}{2}}+\int_{\frac{t}{2}}^t\Big)\int_{\mathbb R}\partial_x^kG(x,t;y,s)g_{3s}dyds=:I^{0,k}_{3,1}+I^{0,k}_{3,2}.
\end{align}
Then 
\begin{align}\label{vE:3.73}
I^{0,k}_{3,1}&=  \int_{\mathbb R}\partial_x^kG(x,t;y,\frac{t}{2})g_{3}(y,\frac{t}{2})dy-\int_{\mathbb R}\partial_x^kG(x,t;y,0)g_{3}(y,0)dy-\int_0^{\frac{t}{2}}\int_{\mathbb R} \partial_s\partial_x^kG(x,t;y,s)g_{3}dyds.
\end{align}
It follows from \eqref{vE:2.5} and \eqref{vE:2.10} that
\begin{align}\label{1vE:3.74}
g_3(y,0)=u_1(y)=-\frac{1}{2}y v_1(y)-\frac{1}{2}\int_{-\infty}^{y}v_1dy=-\Big(\frac{y}{2}\int_{-\infty}^{y}v_1(y_1)dy_1\Big)_y=:-\tilde G_y,
\end{align}
which yields that 
\begin{align}\label{vE:3.74}
\int_{\mathbb R}\partial_x^kG(x,t;y,0)g_{3}(y,0)dy=\int_{\mathbb R}\partial_y \partial_x^kG(x,t;y,0)\tilde G dy.
\end{align}
 Substituting \eqref{vE:3.74} into \eqref{vE:3.73}, together with \eqref{1vE:2.17}  and \eqref{vE:3.8}, leads to 
\begin{align}\label{vE:3.75}
\|I^{0,k}_{3,1}\|_{L^2}&= O(1)t^{-\frac{k}{2}}  \| G_D\|_{L^1}\|g_3(\frac{t}{2}) \|_{L^2}+O(1)t^{-\frac{k}{2}-\frac{1}{2}}\| G_D\|_{L^2}\| \tilde G\|_{L^1}+O(1)t^{-1-\frac{k}{2}}\int_0^{\frac{t}{2}}\| G_D\|_{L^2}\|g_3 \|_{L^1}
\notag\\
&= O(1)\delta t^{-\frac{k}{2}-\frac{5}{4}}+O(1)\delta t^{-\frac{k}{2}-\frac{3}{4}}+O(1)\delta t^{-\frac{k}{2}-\frac{5}{4}}\ln (1+t)=O(1)\delta (1+t)^{-\frac{k}{2}-\frac{3}{4}},
\end{align}
where we have used $\|\tilde G\|_{L^1}= O(1)\delta$ due to $|\tilde G| = O(1)\delta e^{-cy^2}$ from \eqref{E:2.13}.
In addition, it follows from \eqref{1vE:2.17} and \eqref{vE:3.20} that
\begin{align}\label{vE:3.76}
\|I^{0,k}_{3,2}\|_{L^2}&=\| \int_{\frac{t}{2}}^t \int_{\mathbb R}G\partial_s\partial_y^{k}g_{3}dyds+\sum\limits_{\beta<{k}}C_{\beta}\int_{\frac{t}{2}}^t \int_{\mathbb R}G\tilde R^1_{{k}-\beta-1}\partial_s\partial_y^{\beta}g_{3}dyds\|_{L^2}
\notag\\
&=O(1)\int_{\frac{t}{2}}^t \| G_D\|_{L^1}\|\partial_s\partial_y^{k}g_3 \|_{L^2}ds+O(1)\delta \sum\limits_{\beta<{k}}(1+t)^{-\frac{k-\beta}{2}}\int_{\frac{t}{2}}^t \| G_D\|_{L^1}\|\partial_s\partial_y^{\beta}g_3 \|_{L^2}ds
\notag\\
&=O(1)\delta (1+t)^{-\frac{k}{2}-\frac{5}{4}}.
\end{align}
Substituting \eqref{vE:3.75}-\eqref{vE:3.76} into \eqref{vE:3.72} yields that for $k\leq 5$,
\begin{align}\label{new7}
\|I^{0,k}_{3}\|_{L^2}=O(1)\delta (1+t)^{-\frac{k}{2}-\frac{3}{4}}.
\end{align}
When $l=1, k\leq 4$, we deduce from \eqref{vE:3.20}-\eqref{vE:3.21} that
\begin{align}\label{4.82}
I^{1,k}_3
&=\int_{\mathbb R}\partial_x^kG(x,t;y,\frac{t}{2})g_{3s}(\frac{t}{2})dy+\int_{\mathbb R}\partial_t\partial_x^kG(x,t;y,\frac{t}{2})g_{3}(y,\frac{t}{2})dy-\int_{\mathbb R}\partial_t\partial_x^kG(x,t;y,0)g_{3}(y,0)dy
\notag\\
&\quad
+\int_0^{\frac{t}{2}}\int_{\mathbb R}\partial_s\partial_t\partial_x^kG(x,t;y,s)g_{3}dyds+\int_{\frac{t}{2}}^t\int_{\mathbb R}G\partial_s^2\partial_y^{k}g_{3}dyds
\notag\\
&\quad+\sum\limits_{\beta<{k}}C_{\beta}\int_{\frac{t}{2}}^t\int_{\mathbb R}G\tilde R^1_{{k}-\beta-1}\partial_s^2\partial_y^{\beta}g_{3}dyds+\sum\limits_{\beta\leq k}C_{\beta}\int_{\frac{t}{2}}^t\int_{\mathbb R}G\tilde h^1_{k-\beta}\partial_s\partial_y^{\beta}g_{3}dyds.
\end{align}
It follows from \eqref{1vE:3.74} that 
\begin{align}\label{1vE:3.78}
\int_{\mathbb R}\partial_t\partial_x^kG(x,t;y,0)g_{3}(y,0)dy=\int_{\mathbb R}\partial_t\partial_y\partial_x^kG(x,t;y,0)\tilde Gdy.
\end{align}
Thus using \eqref{1vE:2.17}, \eqref{vE:3.8}, \eqref{vE:3.23}-\eqref{1vE:3.23} and \eqref{1vE:3.78} gives that for $k\leq 4$,
\begin{align}\label{new8}
\|I^{1,k}_{3}\|_{L^2}=&O(1)t^{-\frac{k}{2}}\|G_D\|_{L^1}\|g_{3s}(\frac{t}{2})\|_{L^2}+O(1)t^{-1-\frac{k}{2}}\|G_D\|_{L^1}\|g_{3}(\frac{t}{2})\|_{L^2}+O(1)t^{-\frac{3}{2}-\frac{k}{2}}\|G_D\|_{L^2}\|\tilde G\|_{L^1}
\notag\\
&+O(1)t^{-2-\frac{k}{2}}\int_0^{\frac{t}{2}}\|G_D\|_{L^2}\|g_3\|_{L^1}ds+O(1)\int_{\frac{t}{2}}^t\|G_D\|_{L^1}\|\partial_s^2\partial_y^{k}g_{3}\|_{L^2}ds
\notag\\
&+O(1)\delta \sum\limits_{\beta<{k}}(1+t)^{-\frac{k-\beta}{2}}\int_{\frac{t}{2}}^t\|G_D\|_{L^1}\|\partial_s^2\partial_y^{\beta}g_{3}\|_{L^2}ds
\notag\\
&+O(1)\delta \sum\limits_{\beta\leq {k}}(1+t)^{-1-\frac{k-\beta}{2}}\int_{\frac{t}{2}}^t\|G_D\|_{L^1}\|\partial_s\partial_y^{\beta}g_{3}\|_{L^2}ds=O(1)\delta (1+t)^{-\frac{k}{2}-\frac{7}{4}}.
\end{align}
In summary, we have from \eqref{new7} and \eqref{new8} that
\begin{align}\label{vE:3.78}
\|I^{l,k}_{3}\|_{L^2}=O(1)\delta (1+t)^{-l-\frac{k}{2}-\frac{3}{4}}.
\end{align}
Therefore, we conclude from \eqref{1vE:3.44}, \eqref{1vE:3.75} and \eqref{vE:3.78} that 
\begin{lemma}\label{new9}
It holds that for $l+k\leq 5$, $l\leq 1$,
\begin{align}\label{new}
\sum_{i=1}^3\|I^{l,k}_{i}\|_{L^2}=O(1)(\delta+\delta_1) (1+t)^{-l-\frac{k}{2}-\frac{3}{4}}.
\end{align}
\end{lemma}

Next, we estimate the nonlinear terms $I_i^{l,k}$ for $i=4,5,6.$ First we have
\begin{lemma}\label{new10}
It holds that for $l+k\leq 5$, $l\leq 1$,
\begin{align}\label{vE:3.60}
\|I^{l,k}_4\|_{L^2}=O(1)(N^2(0)+\delta+\epsilon^2)(1+t)^{-l-\frac{k}{2}-\frac{3}{4}}.
\end{align}
\end{lemma}
\begin{proof}
When $l=0, k\leq 4$, we have
\begin{align}\label{vE:3.40}
I^{0,k}_4=\int_0^t\int_{\mathbb R}\partial_x^kG(x,t;y,s)g_{1y}dyds=\Big(\int_0^{\frac{t}{2}}+\int_{\frac{t}{2}}^t\Big)\int_{\mathbb R}\partial_x^kG(x,t;y,s)g_{1y}dyds=:I^{0,k}_{4,1}+I^{0,k}_{4,2}.
\end{align}
It follows from \eqref{vE:2.8} and the a priori assumption $N(T)\leq \epsilon$ that 
\begin{align}\label{vE:3.42}
\int_0^{\frac{t}{2}}\|g_1 \|_{L^1}ds=O(1)\epsilon^2\int_0^{\frac{t}{2}}(1+s)^{-\frac{5}{2}}ds=O(1)\epsilon^2,
\end{align}
which, together with Lemma \ref{vlemma3.1} and the Young inequality, yields that for $t>1$ 
\begin{align}\label{vE:3.41}
\|I^{0,k}_{4,1}\|_{L^2}&=\|\int_0^{\frac{t}{2}}\int_{\mathbb R}\partial_y\partial_x^kG(x,t;y,s)g_{1}dyds\|_{L^2}
= O(1)\int_0^{\frac{t}{2}}t^{-\frac{k}{2}}(1+t)^{-\frac{1}{2}} \| G\ast g_1 \|_{L^2}ds
\notag\\
&= O(1)\int_0^{\frac{t}{2}}t^{-\frac{k}{2}}(1+t)^{-\frac{1}{2}} \| G\|_{L^2} \|g_1 \|_{L^1}ds=O(1)(1+t)^{-\frac{k}{2}-\frac{3}{4}}\int_0^{\frac{t}{2}}\|g_1 \|_{L^1}ds
\notag\\
&=O(1)\epsilon^2(1+t)^{-\frac{k}{2}-\frac{3}{4}}.
\end{align}
In addition, it follows from \eqref{vE:3.20}, \eqref{vE:3.23} and \eqref{11vE:3.23} that
\begin{align}\label{4.51}
\|I^{0,k}_{4,2}\|_{L^2}&=\|\int_{\frac{t}{2}}^t\int_{\mathbb R}\partial_x^kG(x,t;y,s)g_{1y}dyds\|_{L^2}
\notag\\
&=\|\int_{\frac{t}{2}}^t\int_{\mathbb R}\Big\{(-1)^{k-1}\partial_y^{k-1}\partial_xG+\sum\limits_{\beta<{k-1}}C_{\beta}\partial_y^{\beta}[-\partial_y(G\tilde R^1_{{k}-\beta-2})
+G\tilde R^1_{{k}-\beta-1}]\Big\}g_{1y}dyds\|_{L^2}
\notag\\
&\leq\int_{\frac{t}{2}}^t \| \partial_xG \|_{L^1}\|\partial_y^k g_1 \|_{L^2}ds+\sum\limits_{\beta<{k-1}}O(1)\delta \int_{\frac{t}{2}}^t (1+t)^{-\frac{k-\beta-1}{2}}\| G_D\|_{L^2}\|\partial_y^{\beta+2} g_1 \|_{L^2}ds
\notag\\
&\quad+\sum\limits_{\beta<{k-1}}O(1)\delta \int_{\frac{t}{2}}^t (1+t)^{-\frac{k-\beta}{2}}\| G_D\|_{L^2}\| \partial_y^{\beta+1} g_1 \|_{L^2}ds.
\end{align}
From \eqref{vE:2.8}, $g_1=O(1)V_x^2$. Since $V \in C([0,T);H^{5}(\mathbb R))$, the derivative order of $g_1$ is at most $4$. This is why we consider $k\leq 4$ in \eqref{4.51}.
From  \eqref{2.20} it holds that 
\begin{align}\label{vE:3.44}
\|\partial_y^ng_1(s)\|_{L^2}
&=O(1)(N^2(0)+\delta)(1+s)^{-\frac{n}{2}-\frac{5}{4}},
\end{align}
which yields that 
\begin{align}\label{vE:3.43}
\|I^{0,k}_{4,2}\|_{L^2}&= O(1)(N^2(0)+\delta)(1+t)^{-\frac{k}{2}-\frac{3}{4}}.
\end{align}
Substituting \eqref{vE:3.41} and \eqref{vE:3.43} into \eqref{vE:3.40} gives that for $k\leq 4$,
\begin{align}\label{new1}
\|I^{0,k}_4\|_{L^2}=O(1)[N^2(0)+\delta+\epsilon^2](1+t)^{-\frac{k}{2}-\frac{3}{4}}.
\end{align}
 When $l=1, k\leq 3$, we deduce from \eqref{vE:3.21} that
\begin{align}\label{vE:3.47}
I^{1,k}_4&=\int_{\mathbb R}\partial_x^kG(x,t;y,t)g_{1y}(t)dy+\int_0^t\int_{\mathbb R}\partial_t\partial_x^kG(x,t;y,s)g_{1y}dyds
\notag\\
&=\int_{\mathbb R}\partial_x^kG(x,t;y,\frac{t}{2})g_{1y}(\frac{t}{2})dy+\int_0^{\frac{t}{2}}\int_{\mathbb R}\partial_t\partial_x^kG(x,t;y,s)g_{1y}dyds
\notag\\
&\quad+\int_{\frac{t}{2}}^t\int_{\mathbb R}\partial_x^{k}G\partial_sg_{1y}dyds+\sum\limits_{\beta\leq k}C_{\beta}\int_{\frac{t}{2}}^t\int_{\mathbb R}G\tilde h^1_{k-\beta}\partial_y^{\beta+1}g_{1}dyds.
\end{align}
It follows from Lemma \ref{vlemma3.1} and \eqref{vE:3.44} that
\begin{align}\label{vE:3.48}
\|\int_{\mathbb R}\partial_x^kG(x,t;y,\frac{t}{2})g_{1y}(\frac{t}{2})dy\|_{L^2}&=\| \partial_x^kG(x,t;y,\frac{t}{2})\ast g_{1y}(\frac{t}{2})\|_{L^2}=O(1)t^{-\frac{k}{2}}\|G_D\|_{L^1}\|g_{1y}\|_{L^2}
\notag\\
&=O(1)(N^2(0)+\delta)(1+t)^{-\frac{k}{2}-\frac{7}{4}}
\end{align}
and from \eqref{vE:3.42} that
\begin{align}\label{vE:3.49}
&\|\int_0^{\frac{t}{2}}\int_{\mathbb R}\partial_t\partial_x^kG(x,t;y,s)g_{1y}dyds\|_{L^2}=\|\int_0^{\frac{t}{2}}\int_{\mathbb R}\partial_t\partial_y\partial_x^kG(x,t;y,s)g_{1}dyds\|_{L^2}
\notag\\
=&O(1)t^{-\frac{k}{2}-\frac{3}{2}}\int_0^{\frac{t}{2}}\|G_D \ast g_1\|_{L^2}ds
=O(1)\epsilon^2(1+t)^{-\frac{k}{2}-\frac{7}{4}}.
\end{align}
Similar to \eqref{4.51}, we use \eqref{vE:3.20} and \eqref{11vE:3.23} to get
\begin{align*}
&\int_{\frac{t}{2}}^t\int_{\mathbb R}\partial_x^{k}G\partial_sg_{1y}dyds=\int_{\frac{t}{2}}^t\int_{\mathbb R}\partial_x\Big((-1)^{k-1}\partial_y^{k-1}G+\sum\limits_{\beta<{k-1}}C_{\beta}\partial_y^{\beta}(G\tilde R^1_{{k}-\beta-2})\Big)\partial_sg_{1y}dyds
\notag\\
=&\int_{\frac{t}{2}}^t\int_{\mathbb R}\partial_xG \partial_s\partial_y^{k} g_{1}dyds+\sum\limits_{\beta<{k-1}}C_{\beta}\int_{\frac{t}{2}}^t\int_{\mathbb R}G\tilde R^1_{{k}-\beta-2}\partial_s \partial_y^{\beta+2}g_{1}dyds
\notag\\
&+\sum\limits_{\beta<{k-1}}C_{\beta}\int_{\frac{t}{2}}^t\int_{\mathbb R}G\tilde R^1_{{k}-\beta-1}\partial_s \partial_y^{\beta+1}g_{1}dyds.
\end{align*}
Thus, for the last two terms in the right-hand side of \eqref{vE:3.47}, we use Lemma \ref{vlemma3.1}, \eqref{vE:3.23} and \eqref{vE:3.44} to get
\begin{align}\label{vE:3.50}
&\|\int_{\frac{t}{2}}^t\int_{\mathbb R}\partial_x^{k}G\partial_sg_{1y}dyds+\sum\limits_{\beta\leq k}C_{\beta}\int_{\frac{t}{2}}^t\int_{\mathbb R}G\tilde h^1_{k-\beta}\partial_y^{\beta+1}g_{1}dyds\|_{L^2}
\notag\\
=&\int_{\frac{t}{2}}^t\|\partial_xG\|_{L^1}\|\partial_s\partial_y^kg_1\|_{L^2}ds+O(1)\delta\Big\{\sum\limits_{\beta<{k-1}} (1+t)^{-\frac{k-\beta-1}{2}}\int_{\frac{t}{2}}^t\|G_D\|_{L^1}\|\partial_s\partial_y^{\beta+2}g_1\|_{L^2}ds
\notag\\
&+\sum\limits_{\beta<{k-1}} (1+t)^{-\frac{k-\beta}{2}}\int_{\frac{t}{2}}^t\|G_D\|_{L^1}\|\partial_s\partial_y^{\beta+1}g_1\|_{L^2}ds+\sum\limits_{\beta\leq k} (1+t)^{-1-\frac{k-\beta}{2}}\int_{\frac{t}{2}}^t\|G_D\|_{L^1}\|\partial_y^{\beta+1}g_1\|_{L^2}ds\Big\}
\notag\\
=&O(1)(N^2(0)+\delta)(1+t)^{-\frac{k}{2}-\frac{7}{4}},
\end{align}
where we have used the fact that
\begin{align}\label{vE:3.52}
\|\partial_s\partial_y^ng_1\|_{L^2}
&=O(1)(N^2(0)+\delta)(1+s)^{-\frac{n}{2}-\frac{9}{4}}.
\end{align}
Substituting \eqref{vE:3.48}-\eqref{vE:3.50} into \eqref{vE:3.47} yields that for $k\leq 3$,
\begin{align}\label{new2}
\|I^{1,k}_4\|_{L^2}=O(1)(N^2(0)+\delta+\epsilon^2)(1+t)^{-\frac{k}{2}-\frac{7}{4}}.
\end{align}
When $l=0, k=5$, we have
\begin{align}\label{vE:3.54}
I^{0,5}_4=\int_0^t\int_{\mathbb R}\partial_x^5G(x,t;y,s)g_{1y}dyds=\Big(\int_0^{\frac{t}{2}}+\int_{\frac{t}{2}}^t\Big)\int_{\mathbb R}\partial_x^5G(x,t;y,s)g_{1y}dyds=:I^{0,5}_{4,1}+I^{0,5}_{4,2}.
\end{align}
Similar to \eqref{vE:3.41}, we get
\begin{align}\label{1vE:3.53}
\|I^{0,5}_{4,1}\|_{L^2}= O(1)\epsilon^2(1+t)^{-\frac{5}{2}-\frac{3}{4}}.
\end{align}
Moreover, it follows from \eqref{vE:3.20} that
\begin{align}\label{vE:3.56}
I^{0,5}_{4,2}&=\int_{\frac{t}{2}}^t\int_{\mathbb R}((-1)^{5}\partial_y^{5}G+\sum\limits_{\beta<{5}}C_{\beta}\partial_y^{\beta}(G\tilde R^1_{{4}-\beta}))g_{1y}dyds
\notag\\
&=\int_{\frac{t}{2}}^t\int_{\mathbb R}\partial_y^{2}G\partial_y^4g_1dyds+\sum\limits_{\beta<{5}}C_{\beta}\int_{\frac{t}{2}}^t\int_{\mathbb R}(\partial_yG\tilde R^1_{{4}-\beta}+G\partial_y\tilde R^1_{{4}-\beta})\partial_y^{\beta}g_{1}dyds.
\end{align}
The key point to deal with $I^{0,5}_{4,2}$ is to control the $\int_{\frac{t}{2}}^t\int_{\mathbb R}\partial_y^{2}G\partial_y^4g_1dyds$. 
Recall that
\begin{equation}\label{vE:3.57}
\partial_y^2G(x,t;y,s)=-\frac{1}{a} (R_{G}+a_yG_y) +\frac{1}{a}G_s,
\end{equation}
which, together with \eqref{vE:3.9}, yields that
\begin{align}\label{vE:3.58}
\int_{\frac{t}{2}}^t\int_{\mathbb R}\partial_y^{2}G\partial_y^4g_1dyds
=&-\int_{\frac{t}{2}}^t\int_{\mathbb R}\frac{1}{a} (R_{G}+a_yG_y)\partial_y^4g_1dyds+\int_{\frac{t}{2}}^t\int_{\mathbb R}\frac{1}{a}G_s\partial_y^4g_1dyds
\notag\\
=&O(1)\delta \int_{\frac{t}{2}}^t\int_{\mathbb R}\Theta(t,s)E(y,s)G_D(x-y, t-s)\partial_y^4g_1dyds+\int_{\mathbb R}\frac{1}{a}G \partial_y^4g_1dy\Big|_{s=\frac{t}{2}}^{s=t}
\notag\\
&+\int_{\frac{t}{2}}^t\int_{\mathbb R}\Big[\Big(\frac{a_s}{a^2} G-\frac{a_y}{a} G_y\Big)\partial_y^4g_1+\Big(\frac{1}{a}G_y-\frac{a_y}{a^2}G\Big)\partial_s\partial_y^3g_1\Big]dyds.
\end{align}
Thus, substituting \eqref{vE:3.58} into \eqref{vE:3.56} and using Lemma \ref{vlemma3.1}, \eqref{vE:3.23}, \eqref{vE:3.44} and \eqref{vE:3.52} lead to 
\begin{align}\label{vE:3.59}
\|I^{0,5}_{4,2}\|_{L^2}=&O(1)\delta \int_{\frac{t}{2}}^t((1+s)^{-1}+(t-s)^{-\frac{1}{2}}(1+s)^{-\frac{1}{2}})\|G_D\|_{L^1}\|\partial_y^4g_1\|_{L^2}ds+O(1)\|\partial_y^{4}g_1(t)\|_{L^2}
\notag\\
&+O(1)\|G_D(\frac{t}{2})\|_{L^1}\|\partial_y^{4}g_1(\frac{t}{2})\|_{L^2}+O(1) \int_{\frac{t}{2}}^t((1+s)^{-\frac{1}{2}}+(t-s)^{-\frac{1}{2}})\|G_D\|_{L^1}\|\partial_s\partial_y^3g_1\|_{L^2}ds
\notag\\
&+O(1)\delta \sum\limits_{\beta<{5}} \int_{\frac{t}{2}}^t((1+s)^{-\frac{1}{2}}+(t-s)^{-\frac{1}{2}})(1+t)^{-\frac{5-\beta}{2}}\|G_D\|_{L^1}\|\partial_y^{\beta}g_1\|_{L^2}ds
\notag\\
=&O(1)(N^2(0)+\delta)(1+t)^{-\frac{5}{2}-\frac{3}{4}},
\end{align}
where we have used the fact from \eqref{vE:3.3} that 
\begin{align*}
\|\int_{\mathbb R}\frac{1}{a}G(x,t;y,t) \partial_y^4g_1(y,t)dy\|_{L^2}=\|\frac{1}{a} \partial_y^4g_1(x,t)\|_{L^2}.
\end{align*}
Substituting \eqref{1vE:3.53} and \eqref{vE:3.59} into \eqref{vE:3.54} yields that 
\begin{align}\label{new3}
\|I^{0,5}_4\|_{L^2}=O(1)(N^2(0)+\delta+\epsilon^2)(1+t)^{-\frac{5}{2}-\frac{3}{4}}.
\end{align}
For the last case with $l=1, k=4$ for the term $I_4^{l,k}$, we deduce from \eqref{vE:3.21} and \eqref{vE:3.23} that
\begin{align}\label{vE:3.61}
I^{1,4}_4&=\int_{\mathbb R}\partial_x^4G(x,t;y,t)g_{1y}(t)dy+\int_0^t\int_{\mathbb R}\partial_t\partial_x^4G(x,t;y,s)g_{1y}dyds
\notag\\
&=\int_{\mathbb R}\partial_x^4G(x,t;y,\frac{t}{2})g_{1y}(\frac{t}{2})dy+\int_0^{\frac{t}{2}}\int_{\mathbb R}\partial_t\partial_x^4G(x,t;y,s)g_{1y}dyds
\notag\\
&\quad+\int_{\frac{t}{2}}^t\int_{\mathbb R}\partial_x^{4}G\partial_sg_{1y}dyds+\sum\limits_{\beta\leq 4}C_{\beta}\int_{\frac{t}{2}}^t\int_{\mathbb R}(\partial_yG\tilde h^1_{4-\beta}+G\partial_y\tilde h^1_{4-\beta})\partial_y^{\beta}g_{1}dyds.
\end{align}
Similar to \eqref{vE:3.48}-\eqref{vE:3.49}, we get
\begin{align}\label{vE:3.62}
\|\int_{\mathbb R}\partial_x^4G(x,t;y,\frac{t}{2})g_{1y}(\frac{t}{2})dy\|_{L^2}=O(1)(N^2(0)+\delta)(1+t)^{-2-\frac{7}{4}}
\end{align}
and 
\begin{align*}
\|\int_0^{\frac{t}{2}}\int_{\mathbb R}\partial_t\partial_x^4G(x,t;y,s)g_{1y}dyds\|_{L^2}
=O(1)\epsilon^2(1+t)^{-2-\frac{7}{4}}.
\end{align*}
Using \eqref{vE:3.57} we have
\begin{align*}
\int_{\frac{t}{2}}^t\int_{\mathbb R}\partial_x^{4}G\partial_sg_{1y}dyds=&\int_{\frac{t}{2}}^t\int_{\mathbb R}\partial_y^{2}G\partial_s \partial_y^3g_1dyds+\sum\limits_{\beta<{4}}C_{\beta}\int_{\frac{t}{2}}^t\int_{\mathbb R}(\partial_yG\tilde R^1_{{3}-\beta}+G\partial_y\tilde R^1_{{3}-\beta})\partial_s\partial_y^{\beta}g_{1}dyds
\notag\\
=&-\int_{\frac{t}{2}}^t\int_{\mathbb R}\frac{1}{a} (R_{G}+a_yG_y)\partial_s\partial_y^3g_1dyds
+\int_{\mathbb R}\frac{1}{a}G\partial_s \partial_y^3g_1dy\Big|_{s=\frac{t}{2}}^{s=t}
\notag\\
&+\int_{\frac{t}{2}}^t\int_{\mathbb R}\Big(\frac{1}{a^2}a_s G\partial_s \partial_y^3g_1+\frac{1}{a}G_y \partial_s^2\partial_y^2g_1-\frac{1}{a^2}a_yG \partial_s^2\partial_y^2g_1\Big)dyds
\notag\\
&+\sum\limits_{\beta<{4}}C_{\beta}\int_{\frac{t}{2}}^t\int_{\mathbb R}(\partial_yG\tilde R^1_{{3}-\beta}+G\partial_y\tilde R^1_{{3}-\beta})\partial_s\partial_y^{\beta}g_{1}dyds.
\end{align*}
From \eqref{2.20}-\eqref{2.22} and \eqref{vE:3.3} we have
\begin{align*}
\|\int_{\mathbb R}\frac{1}{a}G(x,t;y,t) \partial_s\partial_y^3g_1(y,t)dy\|_{L^2}=\|\frac{1}{a} \partial_s\partial_y^3g_1(x,t)\|_{L^2}
\end{align*}
and
\begin{align*}
\|\partial_s^2\partial_y^ng_1\|_{L^2}
=&O(1)(N^2(0)+\delta)(1+s)^{-\frac{n}{2}-\frac{13}{4}},
\end{align*}
which, together with Lemmas \ref{vlemma3.1}-\ref{vlemma3.2}, \eqref{vE:3.23}-\eqref{1vE:3.23}, \eqref{vE:3.44} and \eqref{vE:3.52}, yields that  for the last two terms of \eqref{vE:3.61},
\begin{align}\label{vE:3.65}
&\|\int_{\frac{t}{2}}^t\int_{\mathbb R}\partial_x^{4}G\partial_sg_{1y}dyds+\sum\limits_{\beta\leq 4}C_{\beta}\int_{\frac{t}{2}}^t\int_{\mathbb R}(\partial_yG\tilde h^1_{4-\beta}+G\partial_y\tilde h^1_{4-\beta})\partial_y^{\beta}g_{1}dyds\|_{L^2}
\notag\\
=&O(1)\delta \int_{\frac{t}{2}}^t((1+s)^{-1}+(t-s)^{-\frac{1}{2}}(1+s)^{-\frac{1}{2}})\|G_D\|_{L^1}\|\partial_s\partial_y^3g_1\|_{L^2}ds+O(1)\|\partial_s\partial_y^3g_1(t)\|_{L^2}
\notag\\
&+O(1)\|G_D(\frac{t}{2})\|_{L^1}\|\partial_s\partial_y^3g_1(\frac{t}{2})\|_{L^2}+O(1) \int_{\frac{t}{2}}^t((1+s)^{-\frac{1}{2}}+(t-s)^{-\frac{1}{2}})\|G_D\|_{L^1}\|\partial_s^2\partial_y^2g_1\|_{L^2}ds
\notag\\
&+O(1)\delta \sum\limits_{\beta<{4}} \int_{\frac{t}{2}}^t((1+s)^{-\frac{1}{2}}+(t-s)^{-\frac{1}{2}})(1+t)^{-\frac{4-\beta}{2}}\|G_D\|_{L^1}\|\partial_s\partial_y^{\beta}g_1\|_{L^2}ds
\notag\\
&+O(1)\delta \sum\limits_{\beta\leq {4}} \int_{\frac{t}{2}}^t((1+s)^{-\frac{1}{2}}+(t-s)^{-\frac{1}{2}})(1+t)^{-1-\frac{4-\beta}{2}}\|G_D\|_{L^1}\|\partial_y^{\beta}g_1\|_{L^2}ds
\notag\\
=&O(1)(N^2(0)+\delta)(1+t)^{-2-\frac{7}{4}}.
\end{align}
Substituting \eqref{vE:3.62}-\eqref{vE:3.65} into \eqref{vE:3.61} leads to
\begin{align}\label{new4}
\|I^{1,4}_4\|_{L^2}=O(1)(N^2(0)+\delta+\epsilon^2)(1+t)^{-2-\frac{7}{4}}.
\end{align}
Thus we obtain \eqref{vE:3.60} from \eqref{new1}, \eqref{new2}, \eqref{new3} and \eqref{new4}. Therefore, the proof is completed.
\end{proof}

Next, for the term $I^{l,k}_5$ we have
\begin{lemma}\label{new11}
It holds that for $l+k\leq 5$, $l\leq 1$,
\begin{align}\label{2vE:3.97}
\|I^{l,k}_5\|_{L^2}
=&O(1)(\sqrt{N^2(0)+\delta}+\delta_1) (1+t)^{-l-\frac{k}{2}-\frac{3}{4}}.
\end{align}
\end{lemma}
\begin{proof}
When $l=0, k\leq 4$, we have
\begin{align}\label{vE:3.85}
I^{0,k}_5&=\Big(\int_0^{\frac{t}{2}}+\int_{\frac{t}{2}}^t\Big)\int_{\mathbb R}\partial_x^kG(x,t;y,s)V_{ss}dyds=:I^{0,k}_{5,1}+I^{0,k}_{5,2}.
\end{align}
It follows that
\begin{align*}
I^{0,k}_{5,1}
&=\int_{\mathbb R}\partial_x^kG(x,t;y,\frac{t}{2})V_s(\frac{t}{2})dy-\int_{\mathbb R}\partial_x^kG(x,t;y,0)V_s(y,0)dy-\int_0^{\frac{t}{2}}\int_{\mathbb R}\partial_s\partial_x^kG(x,t;y,s)V_sdyds.
\end{align*}
Note that $V_s(y,0)=\partial_y \tilde V_s(y,0)$ and $\tilde V_s(y,0)\in L^2\cap L^1$ with $\|\tilde V_s(y,0)\|_{L^1}= O(1)\delta_1$,  thus we have
\begin{align*}
\| \partial_x^k G(x,t;y,0)\ast V_s(y,0)\|_{L^2}=\|\partial_x^k \partial_yG(x,t;y,0)\|_{L^2} \|\tilde V_s(y,0)\|_{L^1}=O(1)\delta_1 (1+t)^{-\frac{k}{2}-\frac{3}{4}},
\end{align*}
 which, together with \eqref{2.21} and \eqref{vE:3.8}, yields that
 \begin{align}\label{vE:3.88}
 \|I^{0,k}_{5,1}\|_{L^2}&=O(1)t^{-\frac{k}{2}}\|{G_D}\|_{L^1}\|{V_s}(\frac{t}{2})\|_{L^2}+O(1)\delta_1t^{-\frac{k}{2}-\frac{3}{4}}+O(1)t^{-\frac{k}{2}-1}\int_0^{\frac{t}{2}}\|{G_D}\|_{L^1}\|{V_s}\|_{L^2}ds
 \notag\\
 &=O(1)(\sqrt{N^2(0)+\delta}+\delta_1) (1+t)^{-\frac{k}{2}-\frac{3}{4}}.
 \end{align}
Moreover, it follows from \eqref{2.22}, \eqref{vE:3.20}, \eqref{vE:3.23} and \eqref{11vE:3.23} that for $k\leq 3$, 
\begin{align}\label{vE:3.89}
\|I^{0,k}_{5,2}\|_{L^2}&=\|\int_{\frac{t}{2}}^t\int_{\mathbb R}((-1)^{k}\partial_y^{k}G+\sum\limits_{\beta<{k}}C_{\beta}\partial_y^{\beta}(G\tilde R^1_{{k}-\beta-1}))V_{ss}dyds\|_{L^2}
\notag\\
&= O(1)\int_{\frac{t}{2}}^t\|G_D\|_{L^1}\|\partial_s^2\partial_y^kV\|_{L^2}ds+\sum\limits_{\beta<{k}}O(1)(1+t)^{-\frac{k-\beta}{2}}\int_{\frac{t}{2}}^t\|G_D\|_{L^1}\|\partial_s^2\partial_y^{\beta}V\|_{L^2}ds
\notag\\
&=O(1)\sqrt{N^2(0)+\delta} (1+t)^{-\frac{k}{2}-1},
\end{align}
and 
\begin{align}\label{vE:3.90}
\|I^{0,4}_{5,2}\|_{L^2}
&=\|\int_{\frac{t}{2}}^t\int_{\mathbb R}-\partial_x\partial_y^{3}G+\sum\limits_{\beta<{3}}C_{\beta}\partial_y^{\beta}[-\partial_y(G\tilde R^1_{{2}-\beta})
+G\tilde R^1_{{3}-\beta}]V_{ss}dyds\|_{L^2}
\notag\\
&= O(1)\int_{\frac{t}{2}}^t\|\partial_xG\|_{L^1}\|\partial_s^2\partial_y^3V\|_{L^2}ds+\sum\limits_{\beta<{3}}O(1)(1+t)^{-\frac{3-\beta}{2}}\int_{\frac{t}{2}}^t\|G_D\|_{L^1}\|\partial_s^2\partial_y^{\beta+1}V\|_{L^2}ds
\notag\\
&\quad +\sum\limits_{\beta<{3}}O(1)(1+t)^{-\frac{3-\beta}{2}-\frac{1}{2}}\int_{\frac{t}{2}}^t\|G_D\|_{L^1}\|\partial_s^2\partial_y^{\beta}V\|_{L^2}ds
=O(1)\sqrt{N^2(0)+\delta} (1+t)^{-3}.
\end{align}
Substituting \eqref{vE:3.88}-\eqref{vE:3.90} into \eqref{vE:3.85} gives that for $k\leq 4$, 
\begin{align}\label{new13}
\|I^{0,k}_{5}\|_{L^2}&=O(1)(\sqrt{N^2(0)+\delta}+\delta_1) (1+t)^{-\frac{k}{2}-\frac{3}{4}}.
\end{align}
When $l=1, k\leq 3$, it follows from \eqref{vE:3.20}-\eqref{vE:3.21} that
\begin{align*}
I^{1,k}_5
=&\int_{\mathbb R}\partial_x^kG(x,t;y,\frac{t}{2})V_{ss}(\frac{t}{2})dy
+\int_{\mathbb R}\partial_t\partial_x^kG(x,t;y,\frac{t}{2})V_{s}(y,\frac{t}{2})dy-\int_{\mathbb R}\partial_t\partial_x^kG(x,t;y,0)V_{s}(y,0)dy
\notag\\
&-\int_0^{\frac{t}{2}}\int_{\mathbb R}\partial_s\partial_t\partial_x^kG(x,t;y,s)V_{s}dyds+\int_{\frac{t}{2}}^t\int_{\mathbb R}G\partial_s^3\partial_y^{k}Vdyds+\sum\limits_{\beta<{k}}C_{\beta}\int_{\frac{t}{2}}^t\int_{\mathbb R}G\tilde R^1_{{k}-\beta-1}\partial_s^3\partial_y^{\beta}Vdyds
\notag\\
&
+\sum\limits_{\beta\leq k}C_{\beta}\int_{\frac{t}{2}}^t\int_{\mathbb R}G\tilde h^1_{k-\beta}\partial_s^2\partial_y^{\beta}Vdyds.
\end{align*}
Since $\|\tilde V_s(y,0)\|_{L^1}= O(1)\delta_1$, we have
\begin{align}\label{3vE:3.90}
\| \partial_t\partial_x^k G(x,t;y,0)\ast V_s(y,0)\|_{L^2}=\|\partial_t\partial_x^k \partial_yG(x,t;y,0)\|_{L^2} \|\tilde V_s(y,0)\|_{L^1}=O(1)\delta_1t^{-\frac{k}{2}-\frac{7}{4}},
\end{align}
which, together with \eqref{2.21}-\eqref{2.22}, \eqref{2.25} and \eqref{vE:3.23}-\eqref{1vE:3.23}, yields that for $k\leq 3$,
\begin{align}\label{new14}
\|I^{1,k}_5\|_{L^2}
&=O(1)t^{-\frac{k}{2}}\| G_D(x,t;y,\frac{t}{2})\|_{L^1}\|V_{ss}(\frac{t}{2})\|_{L^2}+O(1)t^{-1-\frac{k}{2}}\| G_D(x,t;y,\frac{t}{2})\|_{L^1}\|V_{s}(\frac{t}{2})\|_{L^2}
\notag\\
&\quad+O(1)\delta_1t^{-\frac{k}{2}-\frac{7}{4}}+O(1)t^{-2-\frac{k}{2}}\int_0^{\frac{t}{2}} \| G_D\|_{L^1}\|V_{s}\|_{L^2}ds
\notag\\
&\quad+O(1)\int_{\frac{t}{2}}^t\| G_D\|_{L^1}\|\partial_s^3\partial_y^{k}V\|_{L^2}ds
+\sum\limits_{\beta<{k}}O(1)(1+t)^{-\frac{k-\beta}{2}}\int_{\frac{t}{2}}^t\| G_D\|_{L^1}\|\partial_s^3\partial_y^{\beta}V\|_{L^2}ds
\notag\\
&\quad
+\sum\limits_{\beta\leq k}O(1)(1+t)^{-1-\frac{k-\beta}{2}}\int_{\frac{t}{2}}^t\| G_D\|_{L^1}\|\partial_s^2\partial_y^{\beta}V\|_{L^2}ds
=O(1)(\sqrt{N^2(0)+\delta}+\delta_1) (1+t)^{-\frac{k}{2}-\frac{7}{4}}.
\end{align}
When $l=0, k=5$, we have
\begin{align}\label{3vE:3.85}
I^{0,5}_5&=\Big(\int_0^{\frac{t}{2}}+\int_{\frac{t}{2}}^t\Big)\int_{\mathbb R}\partial_x^5G(x,t;y,s)V_{ss}dyds=:I^{0,5}_{5,1}+I^{0,5}_{5,2}.
\end{align}
By the same argument as in \eqref{vE:3.88}, we get
 \begin{align}\label{2vE:3.88}
 \|I^{0,5}_{5,1}\|_{L^2}&=O(1)(\sqrt{N^2(0)+\delta}+\delta_1)  (1+t)^{-\frac{5}{2}-\frac{3}{4}}.
 \end{align}
In addition, it follows from \eqref{vE:3.20} that
\begin{align*}
I^{0,5}_{5,2}&=\int_{\frac{t}{2}}^t\int_{\mathbb R}(-\partial_y^{5}G+\sum\limits_{\beta<{5}}C_{\beta}\partial_y^{\beta}(G\tilde R^1_{{4}-\beta}))V_{ss}dyds
\notag\\
&=\int_{\frac{t}{2}}^t\int_{\mathbb R}\partial_y^{2}G\partial_s^2\partial_y^{3}Vdyds+\sum\limits_{\beta<{4}}C_{\beta}\int_{\frac{t}{2}}^t\int_{\mathbb R}G\tilde R^1_{{4}-\beta} \partial_s^2\partial_y^{\beta}Vdyds
\notag\\
&\quad+\int_{\frac{t}{2}}^t\int_{\mathbb R}(\partial_yG\tilde R^1_{{4}-\beta}+G\partial_y\tilde R^1_{{4}-\beta})\partial_s^2\partial_y^{3}Vdyds.
\end{align*}
Then we use \eqref{vE:3.57} to get
\begin{align*}
\int_{\frac{t}{2}}^t\int_{\mathbb R}\partial_y^{2}G\partial_s^2\partial_y^{3}Vdyds
=&-\int_{\frac{t}{2}}^t\int_{\mathbb R}\frac{1}{a} (R_{G}+a_yG_y)\partial_s^2\partial_y^{3}Vdyds+\int_{\frac{t}{2}}^t\int_{\mathbb R}\frac{1}{a}G_s\partial_s^2\partial_y^{3}Vdyds
\notag\\
=&O(1)\delta \int_{\frac{t}{2}}^t\int_{\mathbb R}\Theta(t,s)E(y,s)G_D(x-y, t-s)\partial_s^2\partial_y^{3}Vdyds+\int_{\mathbb R}\frac{1}{a}G \partial_s^2\partial_y^{3}Vdy\Big|_{s=\frac{t}{2}}^{s=t}
\notag\\
&+\int_{\frac{t}{2}}^t\int_{\mathbb R}\Big(\frac{1}{a^2}a_s G\partial_s^2\partial_y^{3}V+\frac{1}{a}G_y \partial_s^3\partial_y^{2}V-\frac{1}{a^2}a_yG \partial_s^3\partial_y^{2}V\Big)dyds,
\end{align*}
which, together with \eqref{2.22}, \eqref{2.25} and Lemma \ref{vlemma3.1}, indicates that
\begin{align}\label{vE:3.94}
&\|\int_{\frac{t}{2}}^t\int_{\mathbb R}\partial_y^{2}G\partial_s^2\partial_y^{3}Vdyds\|_{L^2}
\notag\\
=&O(1)\delta \int_{\frac{t}{2}}^t((1+s)^{-1}+(t-s)^{-\frac{1}{2}}(1+s)^{-\frac{1}{2}})\|G_D\|_{L^1}\| \partial_s^2\partial_y^{3}V\|_{L^2}ds+O(1)\| \partial_s^2\partial_y^{3}V(t)\|_{L^2}
\notag\\
&+O(1)\|G_D(\frac{t}{2})\|_{L^1}\| \partial_s^2\partial_y^{3}V(\frac{t}{2})\|_{L^2}+O(1) \int_{\frac{t}{2}}^t((1+s)^{-\frac{1}{2}}+(t-s)^{-\frac{1}{2}})\|G_D\|_{L^1}\| \partial_s^3\partial_y^{2}V\|_{L^2}ds
\notag\\
=&O(1)\sqrt{N^2(0)+\delta}(1+t)^{-\frac{7}{2}},
\end{align}
where we have used the fact that
\begin{align*}
\|\int_{\mathbb R}\frac{1}{a}G(x,t;y,t) \partial_s^2\partial_y^{3}V(y,t)dy\|_{L^2}=\|\frac{1}{a} \partial_s^2\partial_y^{3}V(x,t)\|_{L^2}.
\end{align*}
Thus, it follows from \eqref{2.22}, \eqref{vE:3.23}, \eqref{11vE:3.23} and \eqref{vE:3.94} that
\begin{align}\label{vE:3.95}
\|I^{0,5}_{5,2}\|_{L^2}=&O(1)\sqrt{N^2(0)+\delta}(1+t)^{-\frac{7}{2}}+O(1)\delta \sum\limits_{\beta<{4}} \int_{\frac{t}{2}}^t(1+t)^{-\frac{5-\beta}{2}}\|G_D\|_{L^1}\|\partial_s^2\partial_y^{\beta}V\|_{L^2}ds
\notag\\
&+O(1)\delta \int_{\frac{t}{2}}^t((1+s)^{-\frac{1}{2}}+(t-s)^{-\frac{1}{2}})(1+t)^{-\frac{5-\beta}{2}}\|G_D\|_{L^1}\|\partial_s^2\partial_y^{3}V\|_{L^2}ds
\notag\\
=&O(1)\sqrt{N^2(0)+\delta}(1+t)^{-\frac{7}{2}}.
\end{align}
Substituting \eqref{2vE:3.88} and \eqref{vE:3.95} into \eqref{3vE:3.85} leads to 
 \begin{align}\label{new15}
 \|I^{0,5}_{5}\|_{L^2}&=O(1)(\sqrt{N^2(0)+\delta}+\delta_1) (1+t)^{-\frac{5}{2}-\frac{3}{4}}.
 \end{align}
When $l=1, k=4$, similar to \eqref{4.82}, we use \eqref{vE:3.20}-\eqref{vE:3.21} to get
\begin{align}\label{1vE:3.97}
I^{1,4}_5
=&\int_{\mathbb R}\partial_x^4G(x,t;y,\frac{t}{2})V_{ss}(\frac{t}{2})dy
+\int_{\mathbb R}\partial_t\partial_x^4G(x,t;y,\frac{t}{2})V_{s}(y,\frac{t}{2})dy-\int_{\mathbb R}\partial_t\partial_x^4G(x,t;y,0)V_{s}(y,0)dy
\notag\\
&-\int_0^{\frac{t}{2}}\int_{\mathbb R}\partial_s\partial_t\partial_x^4G(x,t;y,s)V_{s}dyds+\int_{\frac{t}{2}}^t\int_{\mathbb R}\partial_y^{2}G\partial_s^3\partial_y^2Vdyds+\sum\limits_{\beta<{3}}C_{\beta}\int_{\frac{t}{2}}^t\int_{\mathbb R} G\tilde R^1_{{3}-\beta}\partial_s^3\partial_y^{2}Vdyds
\notag\\
&+\int_{\frac{t}{2}}^t\int_{\mathbb R}(\partial_yG\tilde R^1_{{3}-\beta}+G\partial_y\tilde R^1_{{3}-\beta})\partial_s^3\partial_y^{\beta}Vdyds
+\sum\limits_{\beta\leq 3}C_{\beta}\int_{\frac{t}{2}}^t\int_{\mathbb R}G\tilde h^1_{4-\beta}\partial_s^2\partial_y^{\beta}Vdyds
\notag\\
&+\int_{\frac{t}{2}}^t\int_{\mathbb R}(\partial_yG\tilde h^1_{4-\beta}+G\partial_y\tilde h^1_{4-\beta})\partial_s^2\partial_y^{\beta}Vdyds.
\end{align}
Using \eqref{vE:3.57} again, we have
\begin{align*}
\int_{\frac{t}{2}}^t\int_{\mathbb R}\partial_y^{2}G\partial_s^3\partial_y^{2}Vdyds
=&-\int_{\frac{t}{2}}^t\int_{\mathbb R}\frac{1}{a} (\alpha R_{G}+a_yG_y)\partial_s^3\partial_y^{2}Vdyds+\int_{\frac{t}{2}}^t\int_{\mathbb R}\frac{1}{a}G_s\partial_s^3\partial_y^{2}Vdyds
\notag\\
=&O(1)\delta \int_{\frac{t}{2}}^t\int_{\mathbb R}\Theta(t,s)E(y,s)G_D(x-y, t-s)\partial_s^3\partial_y^{2}Vdyds+\int_{\mathbb R}\frac{1}{a}G \partial_s^3\partial_y^{2}Vdy\Big|_{s=\frac{t}{2}}^{s=t}
\notag\\
&+\int_{\frac{t}{2}}^t\int_{\mathbb R}\Big(\frac{1}{a^2}a_s G\partial_s^3\partial_y^{2}V+\frac{1}{a}G_y \partial_s^4\partial_yV-\frac{1}{a^2}a_yG \partial_s^4\partial_yV\Big)dyds,
\end{align*}
which, together with \eqref{2.25}-\eqref{2.26} and \eqref{vE:3.8}-\eqref{vE:3.9}, indicate that
\begin{align}\label{vE:3.97}
&\|\int_{\frac{t}{2}}^t\int_{\mathbb R}\partial_y^{2}G\partial_s^3\partial_y^{2}Vdyds\|_{L^2}
\notag\\
=&O(1)\delta \int_{\frac{t}{2}}^t((1+s)^{-1}+(t-s)^{-\frac{1}{2}}(1+s)^{-\frac{1}{2}})\|G_D\|_{L^1}\| \partial_s^3\partial_y^{2}V\|_{L^2}ds+O(1)\| \partial_s^3\partial_y^{2}V(t)\|_{L^2}
\notag\\
&+O(1)\|G_D(\frac{t}{2})\|_{L^1}\| \partial_s^3\partial_y^{2}V(\frac{t}{2})\|_{L^2}+O(1) \int_{\frac{t}{2}}^t((1+s)^{-\frac{1}{2}}+(t-s)^{-\frac{1}{2}})\|G_D\|_{L^1}\| \partial_s^4\partial_yV\|_{L^2}ds
\notag\\
=&O(1)\sqrt{N^2(0)+\delta}(1+t)^{-4},
\end{align}
where we have used the fact that
\begin{align*}
\|\int_{\mathbb R}\frac{1}{a}G(x,t;y,t) \partial_s^3\partial_y^{2}V(y,t)dy\|_{L^2}=\|\frac{1}{a} \partial_s^3\partial_y^{2}V(x,t)\|_{L^2}.
\end{align*}
Substituting \eqref{vE:3.97} into \eqref{1vE:3.97}, together with \eqref{2.21}, \eqref{vE:3.8}, \eqref{vE:3.23}-\eqref{11vE:3.23} and \eqref{3vE:3.90}, yields that
\begin{align}\label{new16}
\|I^{1,4}_5\|_{L^2}
=&O(1)(\sqrt{N^2(0)+\delta}+\delta_1) (1+t)^{-2-\frac{7}{4}}.
\end{align}
Thus we obtain \eqref{2vE:3.97} from \eqref{new13}, \eqref{new14}, \eqref{new15} and \eqref{new16}. Therefore, the proof is completed.
\end{proof}

Finally, for the term $I^{l,k}_6$ we have
\begin{lemma}\label{new17}
It holds that for $l+k\leq 5$, $l\leq 1$,
\begin{align}\label{5vE:3.97}
\|I^{l,k}_6\|_{L^2}
=&O(1)\delta \epsilon (1+t)^{-l-\frac{k}{2}-\frac{3}{4}}.
\end{align}
\end{lemma}
\begin{proof}
When $l=0, k\leq 5$, we have
\begin{align*}
I^{0,k}_6&=\Big(\int_0^{\frac{t}{2}}+\int_{\frac{t}{2}}^t\Big)\int_{\mathbb R}\partial_x^kR_{G}(x,t;y,s)Vdyds=:I^{0,k}_{6,1}+I^{0,k}_{6,2}.
\end{align*}
Let $1\leq q_1<2$ in Lemma \ref{vcorollary3.1}, then we use the a priori assumption \eqref{vE:3.37} to get
\begin{align}\label{n3}
\|I^{0,k}_{6,1}\|_{L^2}\leq &\int_0^{\frac{t}{2}} \| \partial_x^kR_{G}\|_{L^1}\|V\|_{L^2}ds=
O(1)\delta\epsilon\int_0^{\frac{t}{2}}[ t^{-\frac{k+1}{2}-\frac{1}{2q_1}}(1+s)^{-\frac{1}{2}+\frac{1}{2q_1}}
+ t^{-\frac{k+1}{2}-\frac{1}{2}}] (1+s)^{-\frac{3}{4}}ds
\notag\\
=&O(1)\delta \epsilon(1+t)^{-\frac{k}{2}-\frac{3}{4}}. 
\end{align}
It follows from the a priori assumption  \eqref{vE:3.37} that for $l\leq 1$, 
\begin{align*}
\|\partial_s^l\partial_y^nV(s)\|_{L^2}&
=O(1)\epsilon(1+s)^{-l-\frac{n}{2}-\frac{3}{4}},
\end{align*}
which implies that 
\begin{align*}
\int_{\frac{t}{2}}^t\Theta(t,s)\|\partial_s^l\partial_y^nV(s)\|_{L^2}ds
&=O(1)\epsilon (1+t)^{-l-\frac{n}{2}-\frac{3}{4}}.
\end{align*}
Together with \eqref{vE:3.9}, \eqref{vE:3.28} and \eqref{vE:3.31}, we have
\begin{align}\label{n4}
\|I^{0,k}_{6,2}\|_{L^2}
=&\|\int_{\frac{t}{2}}^t\int_{\mathbb R}R_G\partial_y^{k} Vdyds+\sum\limits_{\beta<{k}}C_{\beta}\int_{\frac{t}{2}}^t\int_{\mathbb R}G\tilde R^2_{{k}-\beta-1}\partial_y^{\beta}Vdyds\|_{L^2}
\notag\\
=& O(1)\delta \int_{\frac{t}{2}}^t \Theta(t,s) \|G_D\|_{L^1}\|\partial_y^{k} V\|_{L^2}ds+\sum\limits_{\beta<{k}}O(1)\delta (1+t)^{-\frac{k-\beta}{2}}\int_{\frac{t}{2}}^t \Theta(t,s) \|G_D\|_{L^1}\|\partial_y^{\beta} V\|_{L^2} ds
\notag\\
=&O(1)\delta \epsilon (1+t)^{-\frac{k}{2}-\frac{3}{4}}.
\end{align}
Then \eqref{n3} and \eqref{n4} imply that for $k\leq 5$,
\begin{align}\label{new18}
&\|I^{0,k}_{6}\|_{L^2}= \int_0^{\frac{t}{2}} \| \partial_x^kR_{G}\|_{L^1}\|V\|_{L^2}ds=O(1)\delta \epsilon (1+t)^{-\frac{k}{2}-\frac{3}{4}}.
\end{align}
When $l=1, k\leq 4$, it follows from \eqref{vE:3.28}-\eqref{vE:3.29} that
\begin{align*}
I^{1,k}_ 6
=&\int_{\mathbb R}\partial_x^kR_G(x,t;y,\frac{t}{2})V(\frac{t}{2})dy+\int_0^{\frac{t}{2}}\int_{\mathbb R}\partial_t\partial_x^kR_G(x,t;y,s)Vdyds+\int_{\frac{t}{2}}^t\int_{\mathbb R}R_G\partial_s\partial_y^kVdyds
\notag\\
&+\sum\limits_{\beta<{k}}C_{\beta}\int_{\frac{t}{2}}^t\int_{\mathbb R} G\tilde R^2_{{k}-\beta-1}\partial_s\partial_y^{\beta}Vdyds
+\sum\limits_{\beta\leq k}C_{\beta}\int_{\frac{t}{2}}^t\int_{\mathbb R}G\tilde h^2_{k-\beta}\partial_y^{\beta}Vdyds.
\end{align*}
By the same way as in \eqref{n3}-\eqref{n4}, we get that for $k\leq 4$,
\begin{align}\label{new19}
\|I^{1,k}_ 6\|_{L^2}
=&O(1)\delta t^{-1-\frac{k}{2}}\|G_D\|_{L^1}\|V(\frac{t}{2})\|_{L^2}+O(1)\delta \int_{\frac{t}{2}}^t\Theta(t,s)\|G_D\|_{L^1}\|\partial_s\partial_y^k V\|_{L^2}ds
\notag\\
&+\int_{0}^{\frac{t}{2}}\|\partial_t\partial_x^kR_G\|_{L^1}\|V\|_{L^2}ds+O(1)\delta \sum\limits_{\beta<{k}} (1+t)^{-\frac{k-\beta}{2}} \int_{\frac{t}{2}}^t\Theta(t,s)\|G_D\|_{L^1}\|\partial_s\partial_y^{\beta} V\|_{L^2}ds
\notag\\
&+O(1)\delta \sum\limits_{\beta\leq{k}} (1+t)^{-1-\frac{k-\beta}{2}} \int_{\frac{t}{2}}^t\Theta(t,s)\|G_D\|_{L^1}\|\partial_y^{\beta} V\|_{L^2}ds
\notag\\
=&O(1)\delta  \epsilon (1+t)^{-\frac{k}{2}-\frac{7}{4}}.
\end{align}
Thus we obtain \eqref{5vE:3.97} from \eqref{new18}-\eqref{new19}. Therefore, the proof is completed.
\end{proof}

According to the above discussion, it is easy to show the proof of Theorem \ref{theorem 1.1} as follows.
\begin{proof}[\bf Proof of Theorem \ref{theorem 1.1}]
Combining \eqref{4.43}, \eqref{new}, \eqref{vE:3.60}, \eqref{2vE:3.97} and \eqref{5vE:3.97}, we have 
\begin{align*}
\sup \limits_{0\leq t\leq T}\sum_{l+k\leq 5, l\leq 1}\|\partial_t^l\partial_x^k V(\cdot, t)\|_{L^2}\leq C_1[(\sqrt{N^2(0)+\delta}+\delta_1) +(\epsilon+\delta)\epsilon](1+t)^{-l-\frac{k}{2}-\frac{3}{4}},
\end{align*}
which, together with choosing $\epsilon=2C_1(\sqrt{N^2(0)+\delta}+\delta_1) $ small enough, yields that 
\begin{align*}
\sup \limits_{0\leq t\leq T}\sum_{l+k\leq 5, l\leq 1}(1+t)^{l+\frac{k}{2}+\frac{3}{4}}\|\partial_t^l\partial_x^k V(\cdot, t)\|_{L^2}\leq \frac{2}{3}\epsilon<\epsilon.
\end{align*}
Thus, 
the proof of Theorem \ref{theorem 1.1} is completed.
\end{proof}

\

\noindent
{\bf Acknowledgements\ }
 F. Huang's research is supported in part by the National Key R\&D Program of China 2021YFA1000800 and the
National Natural Science Foundation of China (No. 12288201). X. Wu's research is supported in part by the
National Natural Science Foundation of China (No. 12201649).

     \end{document}